\documentclass{article}
\usepackage{graphicx} 
\usepackage{natbib}
\usepackage[utf8]{inputenc}
\usepackage[portrait, margin=1in]{geometry}
\usepackage{float}
\usepackage{amsmath}
\usepackage{amsfonts}
\usepackage{amsthm}
\usepackage{mathtools}
\usepackage{mathabx}
\usepackage{amssymb}

\usepackage{amsrefs}
\usepackage{esvect}
\usepackage{graphicx}
\usepackage{hyperref}
\usepackage[skip=10pt plus1pt, indent=0pt]{parskip}
\usepackage{soul}
\usepackage[x11names]{xcolor}

\newcommand{\dsc}[1]{\colorbox{blue!30!white}{#1}}

\newcommand{\lvb}[1]{{\color{red}#1}}

\newcommand{\D}{\mathcal{D}}
\newcommand{\aD}{\mathcal{D}^*}
\newcommand{\N}{\mathcal{N}}
\newcommand{\e}{\varepsilon}

\newcommand{\RR}
{\mathbb{R}}	
\renewcommand{\Vec}[1]{\mathrm{vec}(#1)}
\DeclareMathOperator{\tr}{tr}			
\DeclareMathOperator{\Div}{div}		

\usepackage[nodisplayskipstretch]{setspace}
\theoremstyle{definition}

\theoremstyle{plain}

\theoremstyle{plain}

\theoremstyle{definition}

\theoremstyle{definition}
\newtheorem*{remark}{Remark}

\usepackage[capitalise]{cleveref}
\usepackage{authblk}

\title{On a Divergence Penalized Landau-de Gennes Model}

\author[1]{Lia Bronsard\footnote{bronsard@mcmaster.ca }}

\author[2]{Jinqi Chen\footnote{jinqi.chen@tufts.edu}}

\author[3]{L\'{e}a Mazzouza\footnote{lea.mazzouza@gmail.com}}

\author[4]{Daniel McDonald\footnote{danmcdon@sas.upenn.edu}}

\author[5]{Nathan Singh\footnote{nesing@ucdavis.edu}}

\author[6]{Dominik Stantejsky\footnote{dominik.stantejsky@univ-lorraine.fr}}

\author[1]{Lee van Brussel \footnote{vanbrulw@mcmaster.ca}}

\affil[1]{Department of Mathematics and Statistics, McMaster University, Hamilton, ON L8S 4L8 Canada}

\affil[2]{Department of Mathematics, Tufts University, Medford, MA 02155 USA}

\affil[3]{Department of Mathematics of Orsay, Universit\'{e} Paris-Saclay, F-91405 Orsay Cedex, France}

\affil[4]{Department of Mathematics, University of Pennsylvania, Philadelphia, PA 19104-6395 USA}

\affil[5]{Department of Mathematics, University of California, Davis, CA 95616 USA}

\affil[6]{Universit\'e de Lorraine, Institut \'Elie Cartan de Lorraine, UMR 7502 CNRS,  54506 Vand\oe uvre-l\`es-Nancy Cedex, France}

\date{\today}

\begin{document}

\maketitle

\begin{abstract}
We give a brief introduction to a divergence penalized Landau-de Gennes functional as a toy model for the study of nematic liquid crystal with colloid inclusion, in the case of unequal elastic constants. We assume that the nematic occupies the exterior of the unit ball, satisfies homeotropic anchoring at the surface of the colloid and approaches a uniform uniaxial state as $|x|\to\infty$. 
We study the 
``small particle'' limit 
and obtain  a representation formula for solutions to the associated Euler-Lagrange equations. 
We also present a numerical analysis of these equations based on a finite element approach and discuss the effect of the divergence penalization on 
the ``Saturn ring'' defects and on 
the properties of the $Q$-tensor. \\
\linebreak
\textbf{Keywords:} Nematic liquid crystal, $Q-$tensors, representation formula, finite elements. \hfill\hfill
\linebreak
\end{abstract}

\section{Introduction}

Loosely speaking, liquid crystals are a material substance which exhibit properties of both isotopic liquid and solid crystalline structures. 
For the purposes of this work, we focus on three-dimensonal samples of \emph{nematic} liquid crystals; a material typically consisting of elongated, rod-like, axially symmetric molecules that maintain long-range orientational order throughout the sample. 
When colloidal particles are immersed into a nematic system, they disturb the orientational order and create topological defects, leading to  remarkable self-assembly phenomena \cite{poulinstarklubenskyweitz97,musevicetal06}, with many potential applications \cite{senyuketal13,porentaetal14}.

In this article we study a nematic surrounding one spherical particle, with homeotropic (i.e.\ the molecules align perpedicular to the particle surface) anchoring at the surface of the colloid, and uniform alignment far away from it. 
We use the tensorial Landau-de Gennes model described below, which has been used by physicists and mathematicians to better describe nematics in settings involving non-orientability, biaxiality, and the presence of line defects (see e.g. ~\cite{ABGL,ACS2021,ACS2024,BaZar,canevari3d,CO1,CO2,DiMiPi,GWZZ,MaZar}). 
This model can be seen in some sense as a relaxation of the non-convex constraints of director models such as the Oseen-Frank director model \cite{BPP,canevari2d,singperturb,OF,MaZar,nguyenzarnescu13}

The Landau-de Gennes model describes liquid crystal through an object called a $Q$-tensor, $Q(x):\Omega\to\mathcal{S}$, where $\Omega\subset\RR^3$ is the region of liquid crystal under consideration and 
$$\mathcal{S}=\{Q\in M_3(\RR):Q_{ij}=Q_{ji},\ \tr(Q)=0\}$$ 
is the set of $3\times 3$, real-valued, symmetric, traceless matrices. 
One of the primary benefits of this matrix representation is that it allows for the distinction between three observed states characteristic of nematic material. 
The first, called the \emph{isotropic} state, occurs when no axis of preferred alignment exists, and in this case, the $Q$-tensor will have all eigenvalues equal to zero. 
The \emph{uniaxial} state is the second state recognized by the $Q$-tensor model. 
In this case, the liquid crystal exhibits a single preferred axis of alignment and is described by $Q$-tensors with two equal, nonzero eigenvalues. 
The preferred axis of alignment is parallel to the normalized principal eigenvector   $m \in \mathbb{S}^2$ and one can write uniaxial $Q$-tensors in the form 
$$Q=s\left(m\otimes m-\frac{1}{3}I\right)$$
where $s\in \RR\setminus\{0\}$ is a constant. 
Lastly, the \emph{biaxial} state is recognized when $Q$ has three distinct eigenvalues. 

The one-constant approximation of the Landau-de Gennes model seeks physically relevant molecular configurations through minimizing the energy $$E_{LdG}(Q)=\int_{\Omega}\left(\frac{L}{2}\left|\nabla Q\right|^2+f_B(Q)\right) dx$$ over $Q$-tensors in an appropriate function space, where $|\cdot|$ is to be interpreted as the Frobenius norm. 
The gradient term serves as an elastic potential with associated elastic constant $L>0$ and the function 
$$f_B(Q)=-\frac a 2\,\tr(Q^2)-\frac b 3\,\tr(Q^3)+\frac c 4\,\tr(Q^2)^2 + d$$ 
is a bulk potential penalizing nonuniaxial $Q$-tensors. 
The constants $a,b,c$ are material-dependent parameters while the constant $d$ is chosen so that $f_B$ is non-negative and attains a global minimum value of zero over the set of uniaxial $Q$-tensors
\begin{equation*}
\left\{s_*\left(v\otimes v-\frac{1}{3}I\right):v\in\mathbb{S}^2\right\}
\end{equation*}
where $s_*=(b+\sqrt{b^2+24ac})/4c>0$.
For simplicity, we choose $s_*=1$ in this article, however, our analysis does not rely on this choice.


Of both physical and mathematical interest is the question of defects. 
Simply put, a defective region in a material sample of nematic liquid crystal is one in which there is a loss of orientational order.
These defects can occur at a single point, or even along one dimensional structures.
It has been shown experimentally that through the inclusion of a spherical colloid with homeotropic anchoring, its radius can dictate the type of defect observed \cite{STARK2001387}.
Indeed, if molecules within the sample are uniformly aligned with respect to some axis far from the inclusion, a large particle radius induces point defects near the colloid, where as a small radius produces a \emph{Saturn ring} defect--a circular ring about the spherical colloid representing an abrupt change in the orientation of the director field. 

The aforementioned defect phenomenon was mathematically verified by Alama, Bronsard \& Lamy in the Landau-de Gennes framework 
\cite{ABL}. 
After appropriate nondimensionalisation, they consider the theoretical setting of nematic material filling all space outside a unit spherical colloid centered at the origin. 
A homeotropic configuration profile is assumed along the surface of the colloid which is stated as the degree one Dirichlet condition $$Q=Q_b:=\left(\nu\otimes \nu-\frac{1}{3}I\right)
\qquad\text{on }\mathbb{S}^2\, ,$$ 
where $\nu\in\mathbb{S}^2$ is the outward unit normal vector along the spherical colloid's surface. 
Away from the particle, $Q$-tensors are to satisfy the degree zero condition 
$$\lim_{|x|\to\infty}Q(x)=Q_{\infty}:=\left(e_z\otimes e_z-\frac{1}{3}I\right)$$ 
which imposes that the preferred axis of alignment coincide with $e_z:=e_3=(0,0,1)$. 
As demonstrated in \cite{ABL}, the Saturn ring behaviour arises in the $Q$-tensor model through an eigenvalue exchange mechanism of the energy minimizing $Q$-tensor. 
This is observed through their derived \emph{explicit} solution for the ``small particle'' limit, which provides a formula for the Saturn ring radius. 
Geometrically, the defect is a consequence of the change in topological structure from the Dirichlet condition $Q_b$ on the colloid's surface to the far-field condition $Q_{\infty}$. 
When the colloid particle is small compared to the characteristic length scale of the nematic liquid crystal, it is the Saturn ring that balances this topological difference.

In this work, we consider unequal elastic constants in the energy and investigate how the elastic constant affects the defect structure for the small particle problem. 
In the case of the large particle problem, see for example \cite{CoKePh}.
Following in the spirit of \cite{ABL}, we wish to derive and analyze a representation formula for solutions of the small particle problem in this novel scenario.
After appropriate non-dimensionalisation (see e.g. \cite{Gartland, ball}), the energy we wish to study is 
\begin{equation*}
{E}(Q)=\int_{\RR^3\setminus B_1(0)}\left(\frac{1}{2}|\nabla Q|^2+\frac{k}{2}|\Div Q|^2+\frac{1}{\xi^2}\widetilde{f_B}(Q)\right)dx
\, ,
\end{equation*}
where $k\in(-1,\infty)$ is related to the elastic constants of the material, and $\xi$ is inversely proportional to the particle radius.
The small particle 
regime corresponds to $\xi\to\infty$, and thus we settle on considering the dimensionless energy
\begin{equation}\label{eq:energy}
E(Q)
\ = \ 
\int_{\Omega}\left(\frac{1}{2}|\nabla Q|^2+\frac{k}{2}|\Div Q|^2\right)dx
\, ,
\end{equation}
where $\Omega=\RR^3\setminus B_1(0)$ and $Q$ is in the affine space $Q\in Q_\infty+H^1(\Omega)$. A $Q$-tensor minimizing the energy $E$ subject to the above homeotropic and far-field conditions will satisfy the Euler-Lagrange equations (see Section~2).
\begin{equation}\label{eq:EL_Qinf}
\begin{cases}
- \Delta Q_{ij}  - \frac{k}{2}\left( \partial_j \Div Q_i + \partial_i \Div Q_j - \frac{2}{3} \text{div}(\Div Q) \delta_{ij}\right) = 0 & \mbox{in}\ \Omega \, ,\\
Q(x)=Q_b\ & \mbox{on}\ \partial\Omega\, ,\\
Q(x)\to Q_{\infty} & \mbox{as}\ |x|\to\infty \, .
\end{cases}
\end{equation}
However, for computational convenience, we replace $Q(x)$ by $Q(x)-Q_{\infty}$ to obtain the equivalent affine problem 
\begin{equation}\label{eq:EL}
\begin{cases}
- \Delta Q_{ij}  - \frac{k}{2}\left( \partial_j \Div Q_i + \partial_i \Div Q_j - \frac{2}{3} \text{div}(\Div Q) \delta_{ij}\right) = 0 & \mbox{in}\ \Omega \, ,\\
Q(x)=Q_b-Q_{\infty}\ & \mbox{on}\ \partial\Omega\, ,\\
|Q(x)|\to 0 & \mbox{as}\ |x|\to\infty \, .
\end{cases}
\end{equation}
so that $Q(x)+Q_{\infty}$ is a solution with our desired boundary behaviour. 

It turns out that treating this system with classic techniques from the field of PDEs such as Fourier transforms and integration against fundamental solutions \cite{Ev} is enough to derive a representation formula for its solution. However, finding  Green's functions for this problem turns out to be analytically challenging and so we finish with a numerical analysis of the system of equations \eqref{eq:EL} using a finite element approach.

\paragraph{Organization of the paper.} In Section \ref{sec:deriv_eqns}, we derive the system of equations. In Section \ref{sec:adjoint}, we derive the adjoint operator associated to the system and find the fundamental solutions to this adjoint system. 
In Section \ref{sec:rep}, we derive the representation formula. 
Finally, in Section \ref{sec:num}, we numerically explore the defect structure of $Q$ for varying elastic parameter $k$.

\paragraph{Acknowledgements:} 
LB is supported by an NSERC (Canada) Discovery Grants.
This work was mostly completed while JC, LM, DM, NS were supported by a Fields Institute Undergraduate Summer Research program.


\section{Derivation of Equations}
\label{sec:deriv_eqns}

The partial differential equations of system \eqref{eq:EL} are the Euler-Lagrange equations for a constrained optimization problem with constraint $\tr(Q)=0$. To obtain the PDE in \eqref{eq:EL}, we consider the modified energy
\begin{equation*}
E_{\lambda}(Q)= \int_{\Omega}\left(\frac{1}{2}|\nabla Q|^2+\frac{k}{2}|\Div Q|^2+\lambda \tr(Q)\right)dx
\end{equation*}
and solve the Lagrange multiplier system
\begin{equation*}
\begin{cases}
\left.\dfrac{d}{dt}\left[E_{\lambda}(Q+t\phi)\right]\right|_{t=0}=0\\[0.5em]
\tr(Q)=0
\end{cases}
\end{equation*}
which must hold for all $\phi\in C_{0}^{\infty}(\Omega;S_3(\RR))$, where $S_3(\RR)$ is the set of $3\times 3$ symmetric matrices with real components.
Indeed, a standard computation of the first variation of $E_{\lambda}$ yields the weak formulation
\begin{equation*}
\sum_{i,j,k=1}^3\int_{\Omega}\partial_kQ_{ij}\partial_k\phi_{ij}\,dx+k\sum_{i=1}^3\int_{\Omega}(\Div Q_i)(\Div \phi_i)\,dx+\int_{\Omega}\lambda \tr(\phi)\,dx=0
\end{equation*}
holding for all $\phi\in C_c^{\infty}(\Omega;S_3(\RR))$. Upon integrating by parts and applying the Fundamental Lemma of the Calculus of Variations, we obtain the pointwise equations 
\begin{equation}\label{eq:ELprelambda}
\begin{cases}
- \Delta Q_{ii}  - k\partial_i \Div Q_i+\lambda = 0 & \mbox{if}\ i=j\\[0.5em]
- \Delta Q_{ij}  - \frac{k}{2}\left( \partial_j \Div Q_i + \partial_i \Div Q_j\right) = 0 & \mbox{if}\ i\neq j
\end{cases}
\end{equation}
holding in $\Omega$. Note that the mixed divergence terms in the second equation arise from the symmetry of $\phi$. To impose the traceless constraint and solve for $\lambda$, we sum along the diagonal and assume $\tr(Q)=0$.
\begin{equation*}
- \sum_{i=1}^3\Delta Q_{ii}  - k\sum_{i=1}^3\partial_i \Div Q_i+3\lambda =-\Delta \tr(Q)  - k\Div\left(\Div Q\right)+3\lambda=- k\Div\left(\Div Q\right)+3\lambda=0
\end{equation*}
giving $\lambda=\frac{k}{3}\Div\left(\Div Q\right)$. Thus, equations \eqref{eq:ELprelambda} with constraint $\tr(Q)=0$ can be written succinctly using the operator notation
\begin{equation*}
\D(Q):=- \Delta Q_{ij}  - \frac{k}{2}\left( \partial_j \Div Q_i + \partial_i \Div Q_j - \frac{2}{3} \text{div}(\Div Q) \delta_{ij}\right) = 0.
\end{equation*}


\section{The Adjoint Operator $\mathcal{D}^{*}$ and Fundamental Solutions}
\label{sec:adjoint}

Consider the punctured domain $\Omega\setminus B_{\e}(p)$, where $B_{\e}(p)$ is a ball of radius $\e>0$ centered at the point $p=(\alpha,\beta,\gamma)\in\Omega$ and suppose $Q$ is a solution of \eqref{eq:EL} in $\Omega\setminus B_{\e}(p)$. 
We observe that our operator $\D$ is not self-adjoint and so to derive an integral representation for $Q$, we first identify the adjoint operator $\aD$.
This can be achieved starting from the equation
\begin{equation}\label{eq:intzero}
0
\ = \ 
\int_{\Omega\setminus B_{\e}(p)} \langle \phi,\mathcal{D}(Q)\rangle\,dx
\, ,
\end{equation}
where $\phi\in C_c^{\infty}({\Omega};S_3(\RR))$ and successively integrating by parts. 
Careful reindexing, i.e.\ sorting terms by the factors $\phi_{ij}$ and $Q_{ij}$, then yields
\begin{align}
\int_{\Omega\setminus B_{\e}(p)}\langle\phi,\D(Q)\rangle\,dx&=\sum_{i,j=1}^3\int_{\Omega\setminus B_{\e}(p)}\left(- \Delta Q_{ij}  - \frac{k}{2}\left( \partial_j \Div Q_i + \partial_i \Div Q_j\right) +\frac{k}{3} \text{div}(\Div Q) \delta_{ij}\right)\phi_{ij}\,dx\nonumber \\
&=\sum_{i,j=1}^3\int_{\Omega\setminus B_{\e}(p)}\left(-\Delta \phi_{ij}-\frac{k}{2}\left(\partial_j \Div \phi_i+\partial_i\Div \phi_j\right)+\frac{k}{3}\partial_i\partial_j\tr(\phi)\right)Q_{ij}\,dx\label{eqn:IPP_adj_w_bdry}\\
&\qquad +\sum_{i,j=1}^3\int_{\Gamma}\left(\frac{\partial\phi_{ij}}{\partial\nu}+\frac{k}{2}\left(\Div(\phi_i)\nu_j+\Div(\phi_j)\nu_i\right)-\frac{k}{3}\partial_j\tr(\phi)\nu_i\right)Q_{ij}\,dS\nonumber\\
&\qquad+\sum_{i,j=1}^3\int_{\Gamma}\left(-\frac{\partial Q_{ij}}{\partial\nu}-\frac{k}{2}\left(\Div(Q_i)\nu_j+\Div(Q_j)\nu_i\right)+\frac{k}{3}\left(\Div(Q)\cdot\nu\right)\delta_{ij}\right)\phi_{ij}\,dS.\nonumber
\end{align}
where $\Gamma = \partial(\Omega\setminus B_{\e}(p))=\partial\Omega\cup \partial B_{\e}(p)$.
Recalling that $\phi$ is compactly supported in $\Omega$, all boundary terms in the above formula vanish.
We end up with the adjoint operator
\begin{equation}\label{eq:adjointeqs}
\mathcal{D}^*(\phi)
\ := \ 
-\Delta \phi_{ij}
-\frac{k}{2}\Big(\partial_j \Div \phi_i+\partial_i\Div \phi_j\Big)
+\frac{k}{3}\partial_i\partial_j\tr(\phi)
\qquad
\text{ in }
\Omega\setminus B_{\e}(p)
\, .
\end{equation}

Next to find the fundamental solutions associated to $\aD$, we proceed by writing our $Q$-tensors in vector form, so that the linear equations \eqref{eq:adjointeqs} can be expressed as matrix multiplication. 
Given a real-valued $3\times 3$ symmetric matrix $Q\in S_3(\RR)$, we associate to $Q$ the vector $\Vec{Q}\in\RR^6$ given by
\begin{equation*}
\Vec{Q}
\ := \
\begin{pmatrix}
Q_{11} \\ Q_{22} \\ Q_{33} \\ Q_{12} \\ Q_{13} \\ Q_{23}
\end{pmatrix}
\qquad
\text{where}
\qquad
Q
\ = \
\begin{pmatrix}
Q_{11} & Q_{12} & Q_{13}\\
Q_{12} & Q_{22} & Q_{23}\\
Q_{13} & Q_{23} & Q_{33}
\end{pmatrix}
\, .
\end{equation*}
Introducing $E^{mn}\in\RR^{3\times 3}$ with $(E^{mn})_{ij}=\delta_{mi}\delta_{nj}$ and $e^{mn}=\Vec{E^{mn}}$, we can write $\Vec{Q} = \sum_{m\leq n=1}^3Q_{mn}e^{mn}$.
Treating $\aD$ as a $6\times 6$ matrix of differential operators, we seek six vectors, denoted by $\Vec{F^{mn}}$, that satisfy the distributional equations
\begin{equation}\label{eq:preFourier}
\aD(\Vec{F^{mn}})=\delta_{(0,0,0)}e^{mn}.
\end{equation}
Upon taking the Fourier transform, equation \eqref{eq:preFourier} has the form $N\, \widehat{\Vec{F^{mn}}}=e^{mn}$ where $N$ is the matrix of spatial frequencies 
\begin{equation*}
N=\widehat{\aD}= \begin{pmatrix}
        |\xi|^2 + \frac{2k}{3} \xi_1^2 
        & -\frac{k}{3} \xi_1^2 
        & - \frac{k}{3} \xi_1^2 &
        k \xi_1 \xi_2 & 
        k \xi_1 \xi_3 & 
        0 \\
        - \frac{k}{3} \xi_2^2 
        & |\xi|^2 + \frac{2k}{3} \xi_2^2 & 
        - \frac{k}{3} \xi_2^2 & 
        k \xi_1 \xi_2 & 
        0 & 
        k \xi_2 \xi_3
        \\
        - \frac{k}{3} \xi_3^2 
        & - \frac{k}{3} \xi_3^2 
        & |\xi|^2 + \frac{2k}{3} \xi_3^2 
        & 0 
        & k \xi_1 \xi_3 
        & k \xi_2 \xi_3 
        \\
        \frac{k}{6} \xi_1 \xi_2 
        & \frac{k}{6} \xi_1 \xi_2 
        & - \frac{k}{3} \xi_1 \xi_2
        & |\xi|^2 + \frac{k}{2} (\xi_1^2 + \xi_2^2) 
        & \frac{k}{2} \xi_2 \xi_3 
        & \frac{k}{2} \xi_1 \xi_3 
        \\
        \frac{k}{6} \xi_1 \xi_3 
        & -\frac{k}{3} \xi_1 \xi_3 & 
        \frac{k}{6} \xi_1 \xi_3 
        & 
        \frac{k}{2} \xi_2 \xi_3 
        & 
        |\xi|^2 + \frac{k}{2} (\xi_1^2 + \xi_3^2) 
        & \frac{k}{2} \xi_1 \xi_2 
        \\
        - \frac{k}{3} \xi_2 \xi_3 & \frac{k}{6} \xi_2 \xi_3
        & \frac{k}{6} \xi_2 \xi_3 
        & \frac{k}{2} \xi_1 \xi_3  
        & \frac{k}{2} \xi_1 \xi_2 
        & |\xi|^2 + \frac{k}{2} (\xi_2^2 + \xi_3^2)
    \end{pmatrix}.
\end{equation*}
Explicit expressions for $\Vec{F^{mn}}$ can be found by solving $\widehat{\Vec{F^{mn}}}=N^{-1}e^{mn}$ and then taking the inverse Fourier transform. 
As an example, $\widehat{\Vec{F^{11}}}$ can be computed to be
\begin{eqnarray*}
    \widehat{\Vec{F^{11}}} &=& \frac{1}{(2 + k) (3 + 2k) (\xi_{1}^{2} + \xi_{2}^{2} + \xi_{3}^{2})^{3}}
    \begin{pmatrix}
    3(2 + k)\xi_{1}^{4} + (12 + k(10 + k))\xi_{1}^{2}(\xi_{2}^{2} + \xi_{3}^{2}) + (2 + k)(3 + 2k) (\xi_{2}^{2} + \xi_{3}^{2})^{2}\\
    k \xi_{2}^{2}(2 ( 1 + k) \xi_{1}^{2} + (2 + k)(\xi_{2}^{2} + \xi_{3}^{2}))\\
    k \xi_{3}^{2}(2 ( 1 + k) \xi_{1}^{2} + (2 + k)(\xi_{2}^{2} + \xi_{3}^{2}))\\
    - k \xi_{1} \xi_{2}(\xi_{1}^{2} + (1 + k)(\xi_{2}^{2} + \xi_{3}^{2}))\\
    - k \xi_{1} \xi_{3} (\xi_{1}^{2} + (1 + k) (\xi_{2}^{2} + \xi_{3}^{2}))\\
    k \xi_{2} \xi_{3} (2 ( 1 + k) \xi_{1}^{2} + (2 + k) (\xi_{2}^{2} + \xi_{3}^{2}))
    \end{pmatrix}
    \, .
    \end{eqnarray*}

A straightforward computation shows that    
\begin{eqnarray*}
\widecheck{\frac{\xi_{1}^3 \xi_{2}}{(\xi_{1}^{2} + \xi_{2}^{2} + \xi_{3}^{2})^{3}}} 
&=& 
\frac{-3 x y (y^2 + z^2)}{32 \pi (x^2 + y^2 + z^2)^{\frac{5}{2}}}\\
\widecheck{\frac{\xi_{1}^{2} \xi_{2} \xi_{3}}{(\xi_{1}^{2} + \xi_{2}^{2} + \xi_{3})^{3}}} 
&=&
\frac{y z (2 x^{2} - y^2 - z^2)}{32 \pi (x^{2} + y^{2} + z^{2})^{\frac{5}{2}}}\\
\widecheck{\frac{\xi_{1}^{4}}{(\xi_{1}^{2} + \xi_{2}^{2} + \xi_{3})^{3}}} 
&=& 
\frac{3 (y^{2} + z^{2})^{2}}{32 \pi (x^{2} + y^{2} + z^{2})^{\frac{5}{2}}}\\
\widecheck{\frac{\xi_{1}^{2} \xi_{2}^{2}}{(\xi_{1}^{2} + \xi_{2}^{2} + \xi_{3})^{3}}} 
&=& 
\frac{3 x^{2} y^{2} + (x^{2} + y^{2} + z^{2}) z^{2} }{32 \pi (x^{2} + y^{2} + z^{2})^{\frac{5}{2}}}
\, ,
\end{eqnarray*}
and therefore $\Vec{F^{11}}$ can be written as
\begin{equation*}
\Vec{F^{11}}
\ = \ 
\frac{1}{32 \pi (3 + 2k) (2 + k)}\left( \frac{48}{r} \begin{pmatrix}
                1 \\
                0 \\
                0 \\
                0 \\
                0 \\
                0 \\
            \end{pmatrix} +
        \frac{1}{r}
        \begin{pmatrix}
                (7k^2+40k)r \\
                (5k^2+8k)r \\
                (5k^2+8k)r \\
                0 \\
                0 \\
                0 \\
            \end{pmatrix} +
        3k^{2}x^{2}
        \begin{pmatrix}
                x^2 \\
                y^2 \\
                z^2 \\
                xy \\
                xz \\
                yz \\
            \end{pmatrix} +
        \frac{1}{r^{3}}
        \begin{pmatrix}
                (6k^2+16k)x^2 \\
                -(k^2 x^2 + (5k^2+8k)y^2) \\
                -(k^2 x^2 + (5k^2+8k)y^2) \\
                k(4+k)xy \\
                k(4+k)xz \\
                -k(8+5k)yz \\
            \end{pmatrix} \right)
\end{equation*}
where $r=\sqrt{x^2+y^2+z^2}$. 
We observe that when $k=0$ we recover the usual fundamental solution for the Laplace operator:
\begin{align*}
\Vec{F^{11}}
\ &= \
\frac{1}{192 \pi} \frac{1}{r^5}
\begin{pmatrix}
48 r^4  \\
0  \\
0  \\
0  \\
0  \\
0
\end{pmatrix}
\ = \
\frac{1}{4 \pi} \frac{1}{r}
\begin{pmatrix}
1  \\
0  \\
0  \\
0  \\
0  \\
0
\end{pmatrix}
\, .
\end{align*}

\section{The Representation Formula}
\label{sec:rep}

This section is dedicated to finding a representation formula for a solution $Q$ of \eqref{eq:EL} in terms of $Q|_{\partial\Omega}$ and the derivatives $(\nabla Q)|_{\partial\Omega}$.

Revisiting the computations leading to equation \eqref{eqn:IPP_adj_w_bdry}, we find it useful to define the operators
\begin{align*}
\mathcal{L}(\phi)&:=\frac{\partial\phi_{ij}}{\partial\nu}+\frac{k}{2}\left(\Div(\phi_i)\nu_j+\Div(\phi_j)\nu_i\right)-\frac{k}{3}\partial_j\tr(\phi)\nu_i\\[0.5em]
\N(Q)&:=-\frac{\partial Q_{ij}}{\partial\nu}-\frac{k}{2}\left(\Div(Q_i)\nu_j+\Div(Q_j)\nu_i\right)+\frac{k}{3}\left(\Div(Q)\cdot\nu\right)\delta_{ij}
\, .
\end{align*}
Inserting these definitions and \eqref{eqn:IPP_adj_w_bdry} into equation \eqref{eq:intzero}, we obtain
\begin{equation}\label{eq:newintzero}
\int_{\Omega\setminus B_{\e}(p)} \langle \aD(\phi),Q\rangle\,dx
+\int_{\Gamma}\langle\mathcal{L}(\phi),Q\rangle\,dS
+\int_{\Gamma}\langle\phi,\N(Q)\rangle\,dS
\ = \ 0
\, .
\end{equation}
Note that if $\phi$ is a solution of $\aD(\phi)=0$ on $\Omega\setminus B_\e(p)$, then the interior integral in equation \eqref{eq:newintzero} vanishes, leaving us only the boundary integrals.
As shown in the previous section, $F^{mn}$ are solutions to $\aD(\phi)=0$ on $\Omega\setminus B_\e(p)$ so we can replace $\phi$ by $F^{mn}_p$ in \eqref{eq:newintzero}, where $F^{mn}_p$ is the $mn$-th fundamental solution for $\aD$ translated to the point $p=(\alpha,\beta,\gamma)\in\Omega$.
By doing so, we are left with the sum
\begin{equation}\label{eq:newnewintzero}
\int_{\partial\Omega}\langle\mathcal{L}(F_p^{mn}),Q_b-Q_{\infty}\rangle\,dS
+\int_{\partial\Omega}\langle F_p^{mn},\N(Q)\rangle\,dS
+\int_{\partial B_{\e}(p)}\langle\mathcal{L}(F_p^{mn}),Q\rangle\,dS
+\int_{\partial B_{\e}(p)}\langle F_p^{mn},\N(Q)\rangle\,dS
\ = \ 
0
\, ,
\end{equation}
where we have used $\Gamma=\partial\Omega\cup \partial B_{\e}(p)$.
We now study the integrals over $\partial B_\e(p)$ and first show that the last term in \ref{eq:newnewintzero} converges to zero as $\e\to 0^+$, i.e.\
\begin{equation*}
\lim_{\e\to0^+}\int_{\partial B_{\e}(p)}\langle F_p^{mn},\N(Q)\rangle\,dS
\ = \ 
0
\, .
\end{equation*}

Indeed, it follows from the calculations of the fundamental solutions $F^{mn}$ that they can be written in the form
\begin{equation*}
F^{mn}_{ij}(x,y,z)
\ = \ 
C^{mn}_{ij}(k)\frac{P_{ij}^{mn}(x,y,z)}{\left(x^2+y^2+z^2\right)^{5/2}}
\, ,
\end{equation*}
where $C^{mn}_{ij}(k)$ are constants depending on the elastic constant $k>0$ and $P_{ij}^{mn}$ denote polynomials of degree four in the spatial variables. 
Moreover, the smoothness of the solution $Q$ near $B_{\e}(p)$ gives $\|\N(Q)\|_{L^{\infty}}<\infty$. 
Thus, there exists a constant $C>0$ independent of $\e$ such that
\begin{equation}\label{eq:secondlastlim}
\left|\int_{\partial B_{\e}(p)}\langle F_p^{mn},\N(Q)\rangle\,dS\right|
\ \leq \ 
\frac{C}{\e}\int_{\partial B_{\e}(p)}\,dS
\ = \ 
4\pi C\e \xrightarrow[\e \to 0^+]{} 0
\end{equation}
as desired.

For the third term in \eqref{eq:newnewintzero}, we remark, since the fundamental solution components are explicitly known, that all $\mathcal{L}(F_p^{mn})_{ij}$ are of the form
\begin{equation*}
\frac{P_{ij}^{mn}(x-\alpha,y-\beta,z-\gamma)}{((x-\alpha)^2+(y-\beta)^2+(z-\gamma)^2)^3}
\, ,
\end{equation*}
where $P$ denotes a polynomial of degree four.
Thus, using the continuity of $Q$, we can write the limit as
\begin{equation*}
\lim_{\e\to 0^+}\int_{\partial B_{\e}(p)}\langle\mathcal{L}(F_p^{mn}),Q\rangle\,dS=\sum_{i,j=1}^3C_{ij}^{mn}Q_{ij}(\alpha,\beta,\gamma)
\, ,
\end{equation*}
where 
\begin{equation*}
C_{ij}^{mn}=\lim_{\e\to0^+}\int_{\partial B_{\e}(p)}\mathcal{L}(F_p^{mn})_{ij}\,dS
\, .
\end{equation*}
Upon integrating in spherical coordinates centered at $p$ with radius $\e$, one can explicitly compute
\begin{equation*}
\int_{\partial B_{\varepsilon}(p)} \frac{(x-\alpha)^{\beta_1} (y-\beta)^{\beta_2} (z-\gamma)^{\beta_3}}{\left[ (x-\alpha)^2 + (y-\beta)^2 + (z-\gamma)^2 \right]^3} \, dS
= \begin{cases}
\displaystyle \dfrac{4\pi}{5}, & \text{if } \exists\, i \text{ such that } \beta_i = 4, \\[0.5em]
\displaystyle \dfrac{4\pi}{15}, & \text{if } \exists\, i \neq j \text{ with } \beta_i = \beta_j = 2, \\[0.5em]
0, & \text{otherwise},
\end{cases}
\end{equation*}
where $\beta_i \in \{0,1,2,3,4\}$ and $\sum_{i=1}^3 \beta_i = 4$.
Hence, using the explicit form of the $F_p^{mn}$, we obtain:
\begin{equation*}
C^{mn}_{ij}=
\begin{cases}
1 & \mbox{if}\ (i,j)=(m,n)\ \mbox{or}\ (n,m)\\
0 & \mbox{otherwise}.
\end{cases}
\end{equation*}
Finally,
\begin{equation}\label{eq:lastlim}
\lim_{\e\to 0^+}\int_{\partial B_{\e}(p)}\langle\mathcal{L}(F_p^{mn}),Q\rangle\,dS
\ = \ 
\lim_{\e\to 0^+}\sum_{i,j=1}^3\int_{\partial B_{\e}(p)}\langle\mathcal{L}(F_p^{mn})_{ij},Q_{ij}\rangle\,dS
\ = \
2Q_{mn}-Q_{mn}\delta_{mn}
\, .
\end{equation}

We conclude from \eqref{eq:newnewintzero} using \eqref{eq:secondlastlim} and \eqref{eq:lastlim} that
\begin{equation*}
(2-\delta_{mn})Q_{mn}(\alpha,\beta,\gamma)
\ = \ 
-\int_{\partial\Omega}\langle F^{mn}_{p},\N(Q)\rangle \ dS-\int_{\partial \Omega}\langle \mathcal{L}(F^{mn}_{p}),Q_b-Q_{\infty}\rangle \ dS
\end{equation*}

\begin{remark}
If one could construct a Green's function $G$ and use it in place of the fundamental solution, the representation formula reduces to
\begin{equation*}
(2-\delta_{mn})Q_{mn}(\alpha,\beta,\gamma)
\ = \
-\int_{\partial \Omega}\langle \mathcal{L}(G^{mn}_{w}),Q_b-Q_{\infty}\rangle \ dS
\, .
\end{equation*}
\end{remark}


\section{Numerical Study of the Euler-Lagrange Equations}\label{sec:num}

In this section, we solve the system of equations \eqref{eq:EL_Qinf} numerically using a finite element approach.
We restrict ourselves to a large box and replace the condition at infinity by a Dirichlet condition on the boundary of this box. 

We generate an adapted mesh with the help of GMSH \cite{gmsh2020} that contains around $130\: 000$ elements. 
Inside the box $[-10,10]^3$ we place a spherical particle of radius $1$ and chose a mesh size $h=0.02$ close to the particle surface to capture the singular structure in detail, while further away, we use $h=0.5$.
We then solve the weak formulation using $P^1-$elements. 
The implementation is carried out in FEniCS \cite{Fenics2015}.
We illustrate our findings graphically by generating images in Paraview \cite{Paraview2015}.

\begin{figure}[H]
\begin{center}
\includegraphics[scale=0.15]{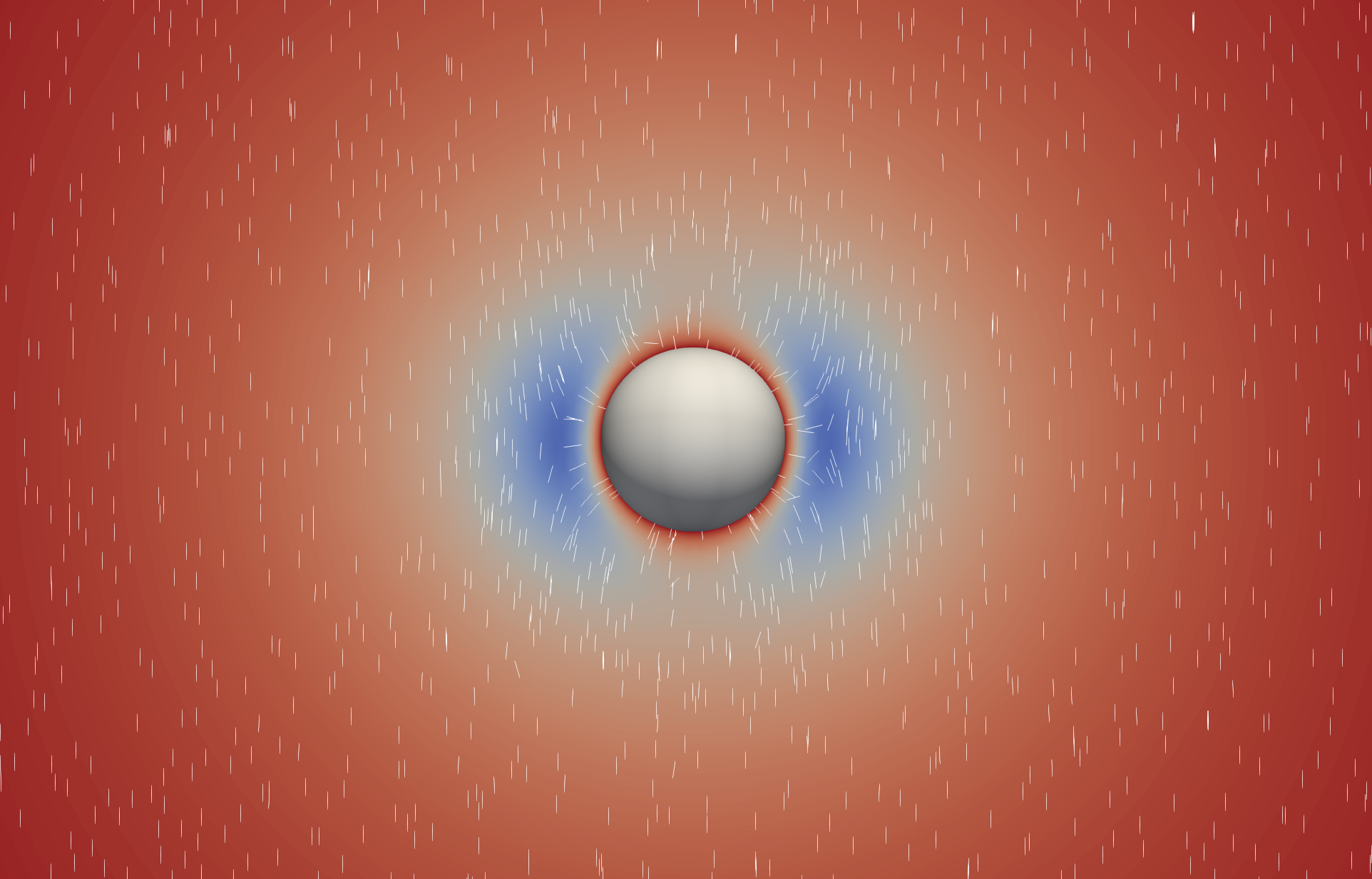} 
\includegraphics[scale=0.15]{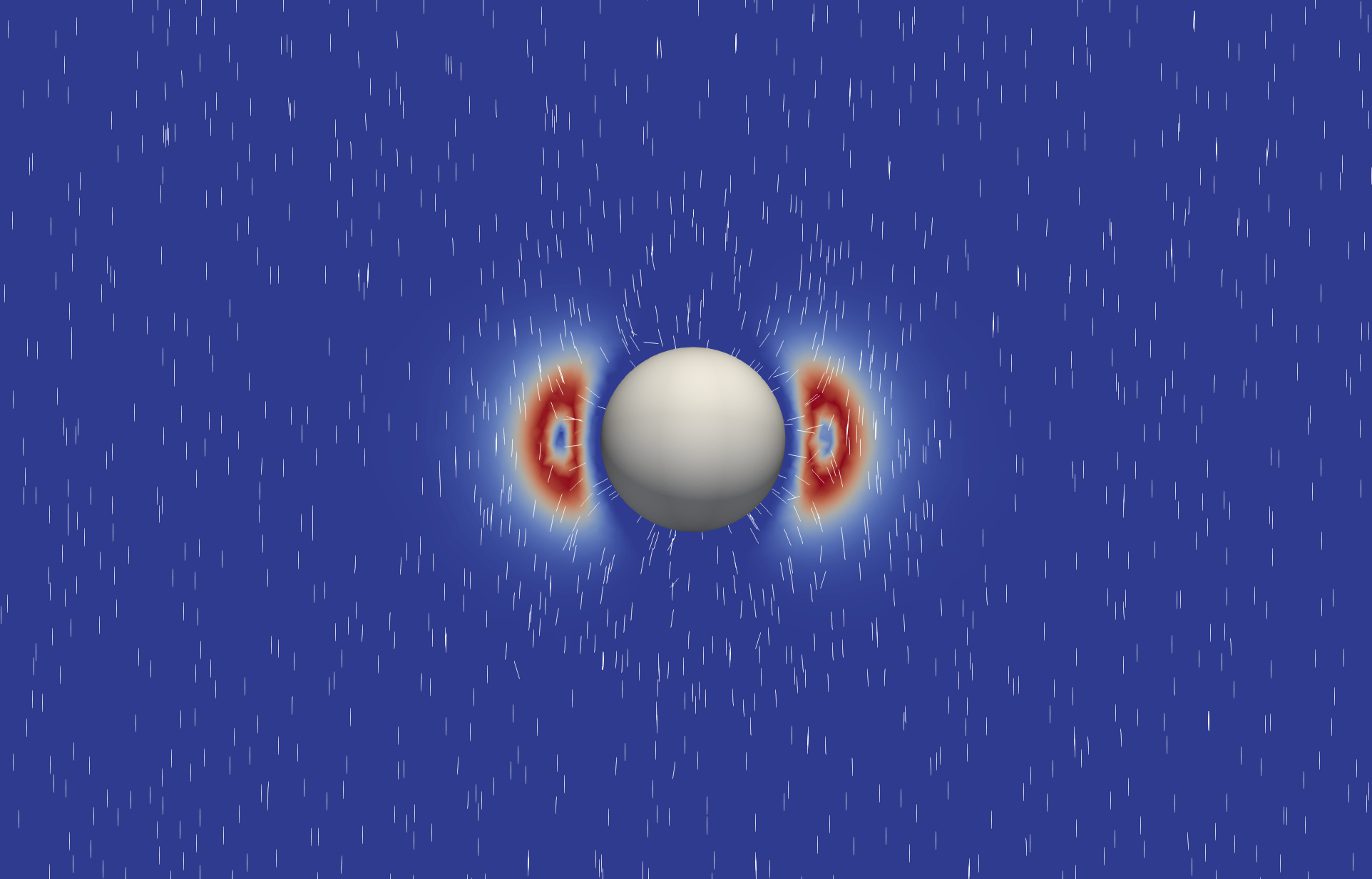} 
\includegraphics[scale=0.15]{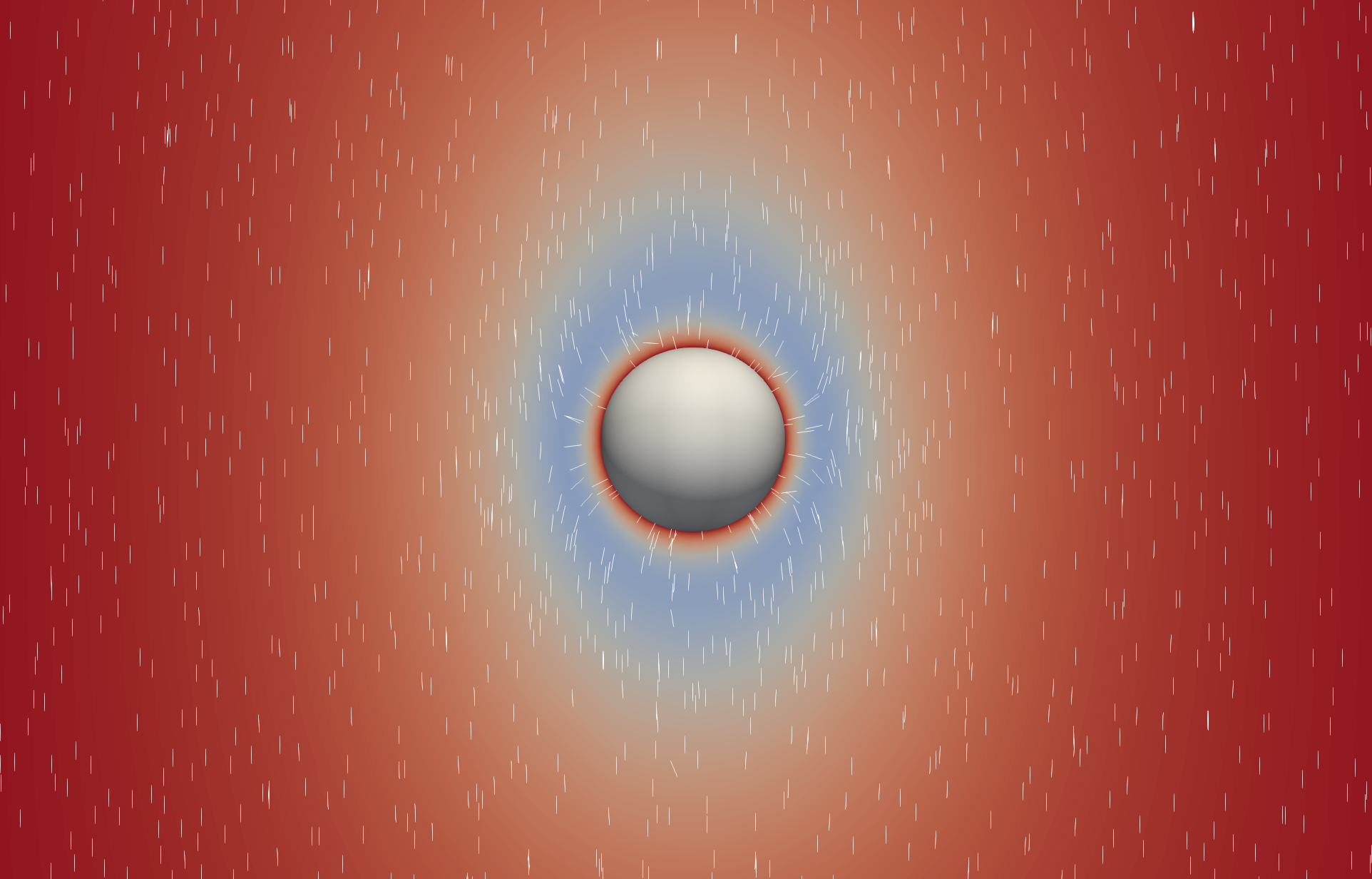} 
\includegraphics[scale=0.15]{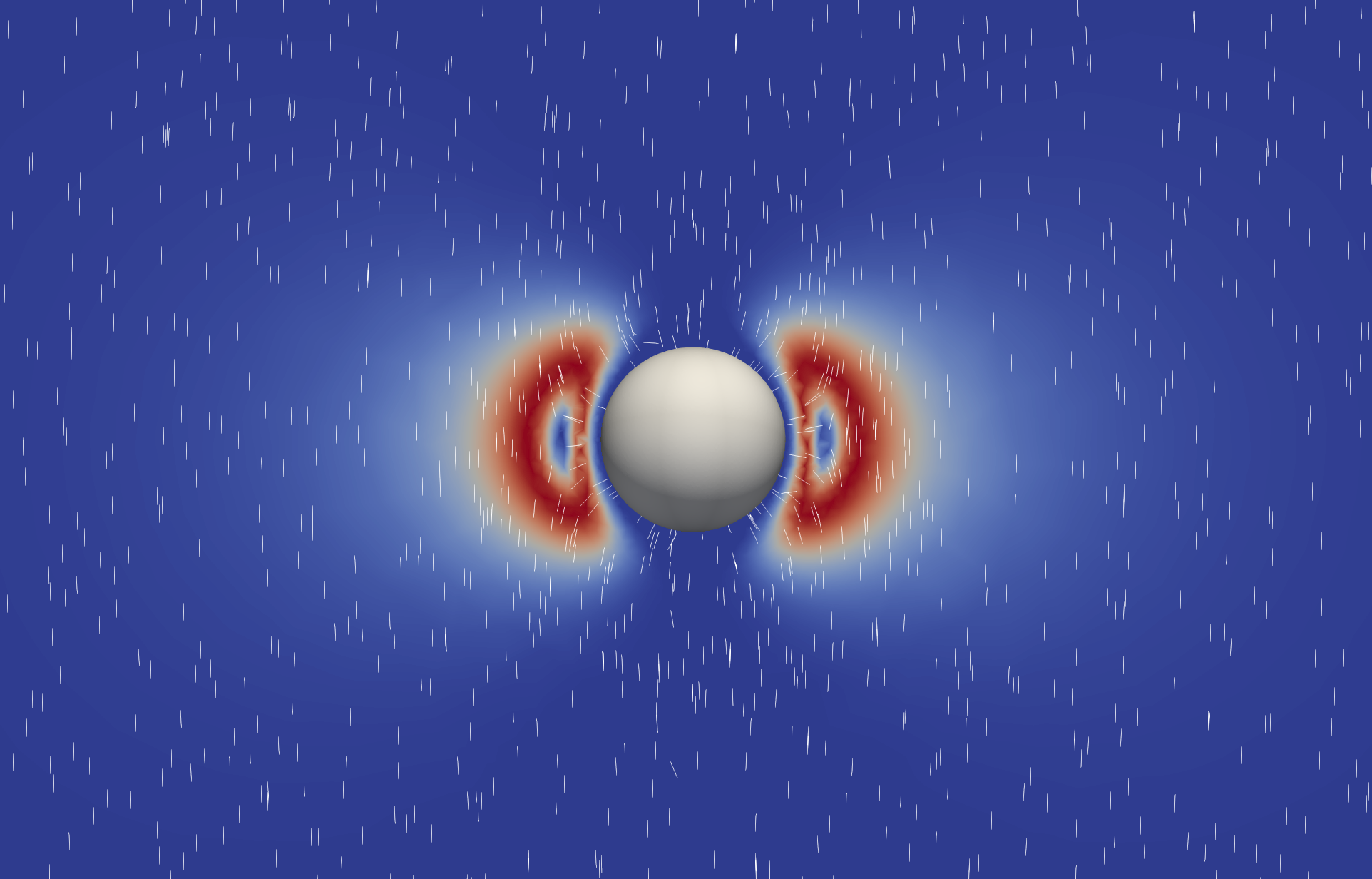} 
\includegraphics[scale=0.15]{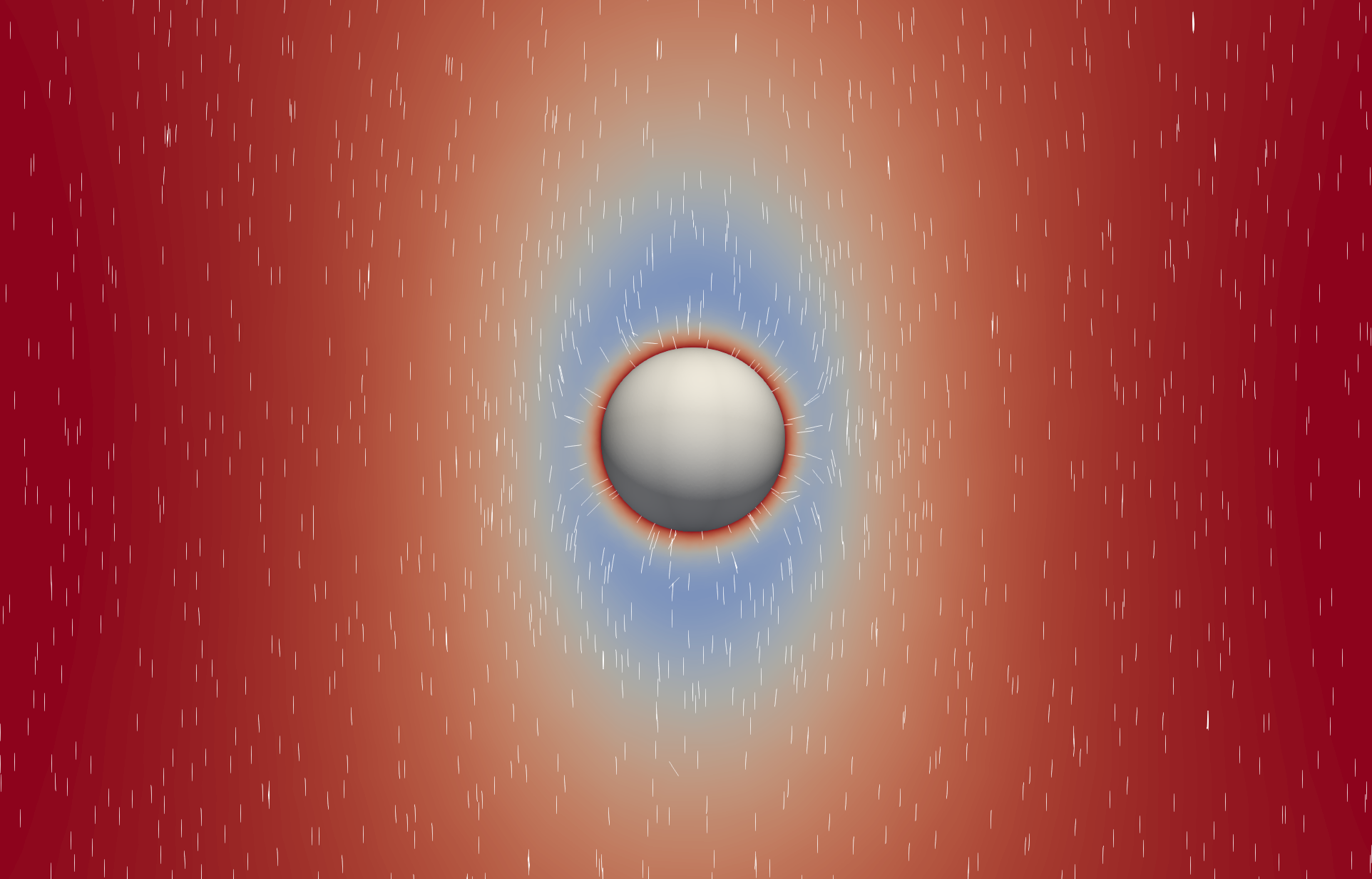} 
\includegraphics[scale=0.15]{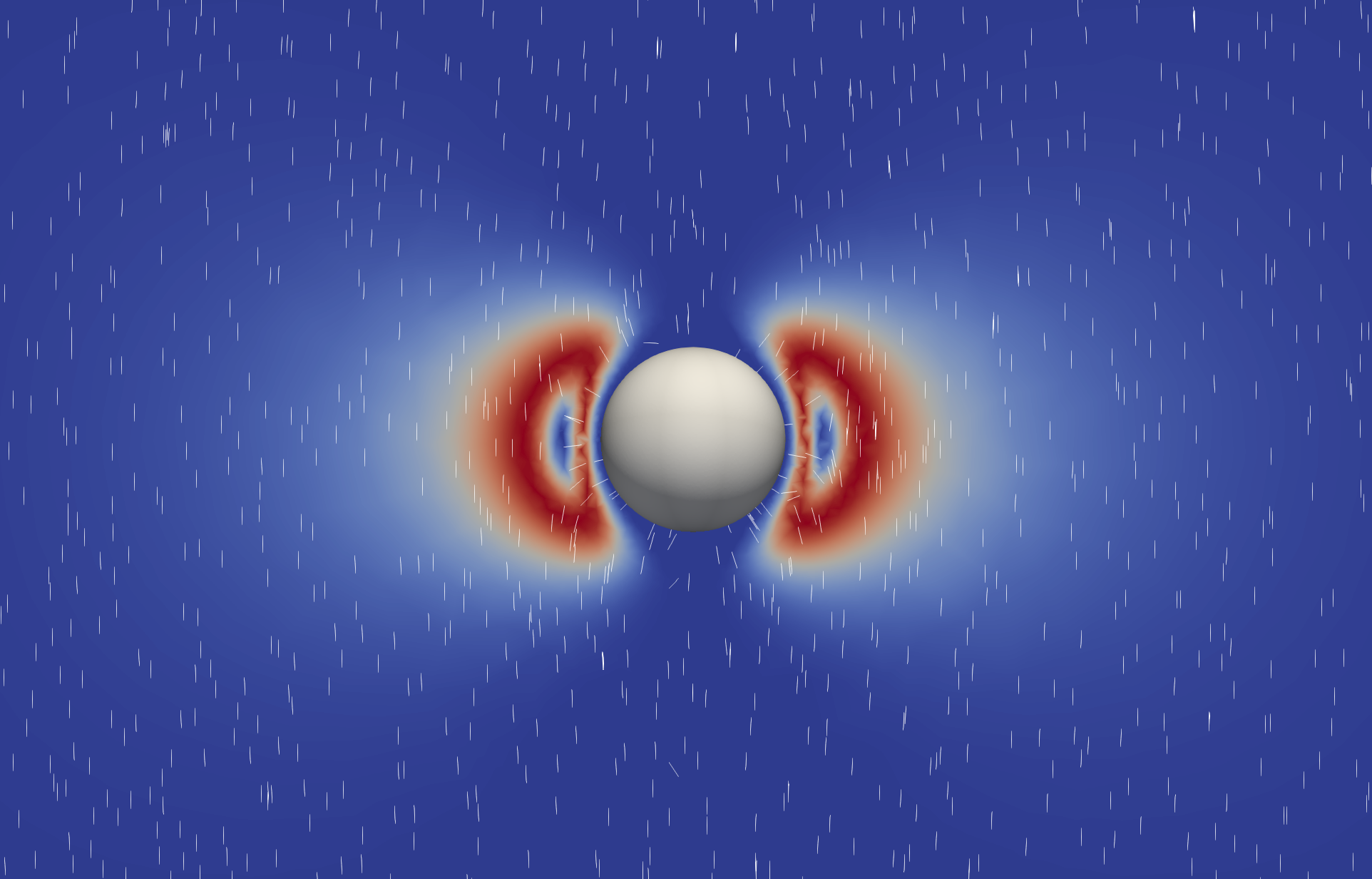} 
\includegraphics[scale=0.15]{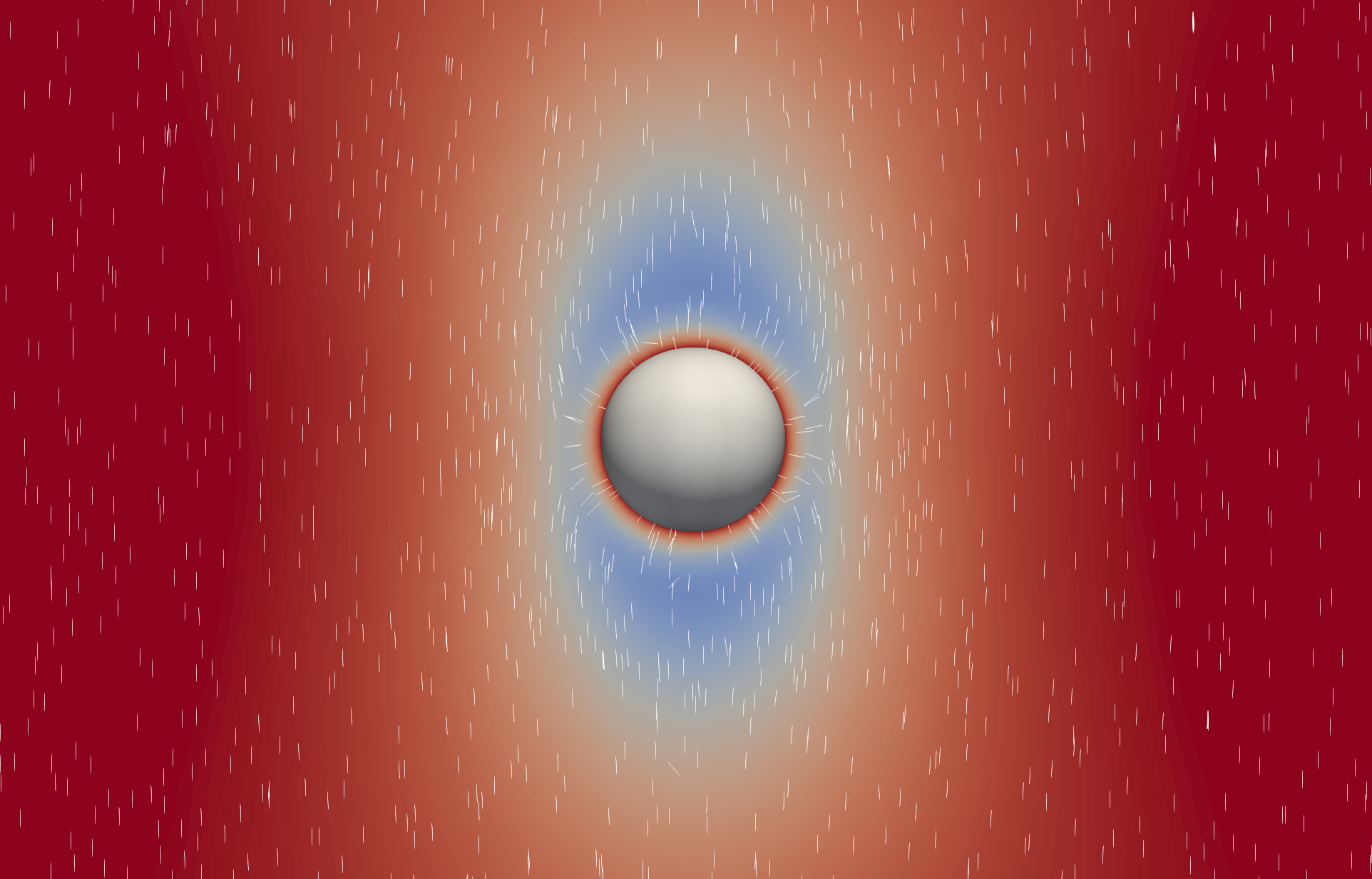} 
\includegraphics[scale=0.15]{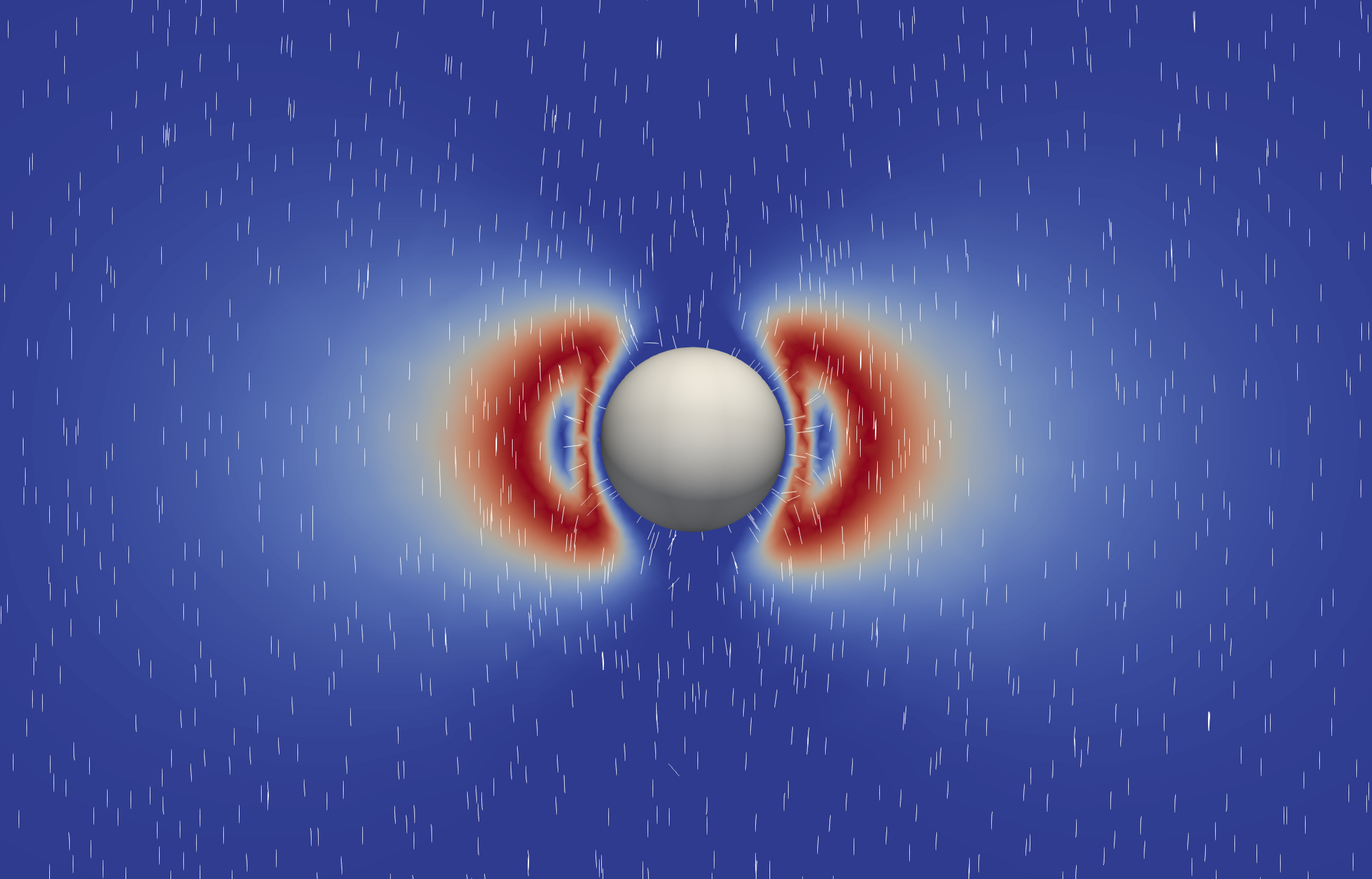} 
\caption{Plots of the norm $|Q|$ (left column) and biaxiality parameter $\beta(Q)$ (right column) for $k=0,5,10,20$ (top to bottom). 
The lines indicate the director field (dominating eigenvector) of $Q$.
The color coding is identical for all four norm images and goes from $0.18$ (dark blue) to $0.82$ (red).
The biaxiality parameter varies between $0$ (uniaxial, dark blue) and $1$ (maximally biaxial, red).
}
\label{fig:norm_biax_k_0_5_10_20}
\end{center}
\end{figure}

\begin{figure}[H]
\begin{center}
\includegraphics[scale=0.15]{Laplace_GradDiv_B1_k0_norm_director.png} 
\includegraphics[scale=0.15]{Laplace_GradDiv_B1_k0_biax_director.png} 
\includegraphics[scale=0.15]{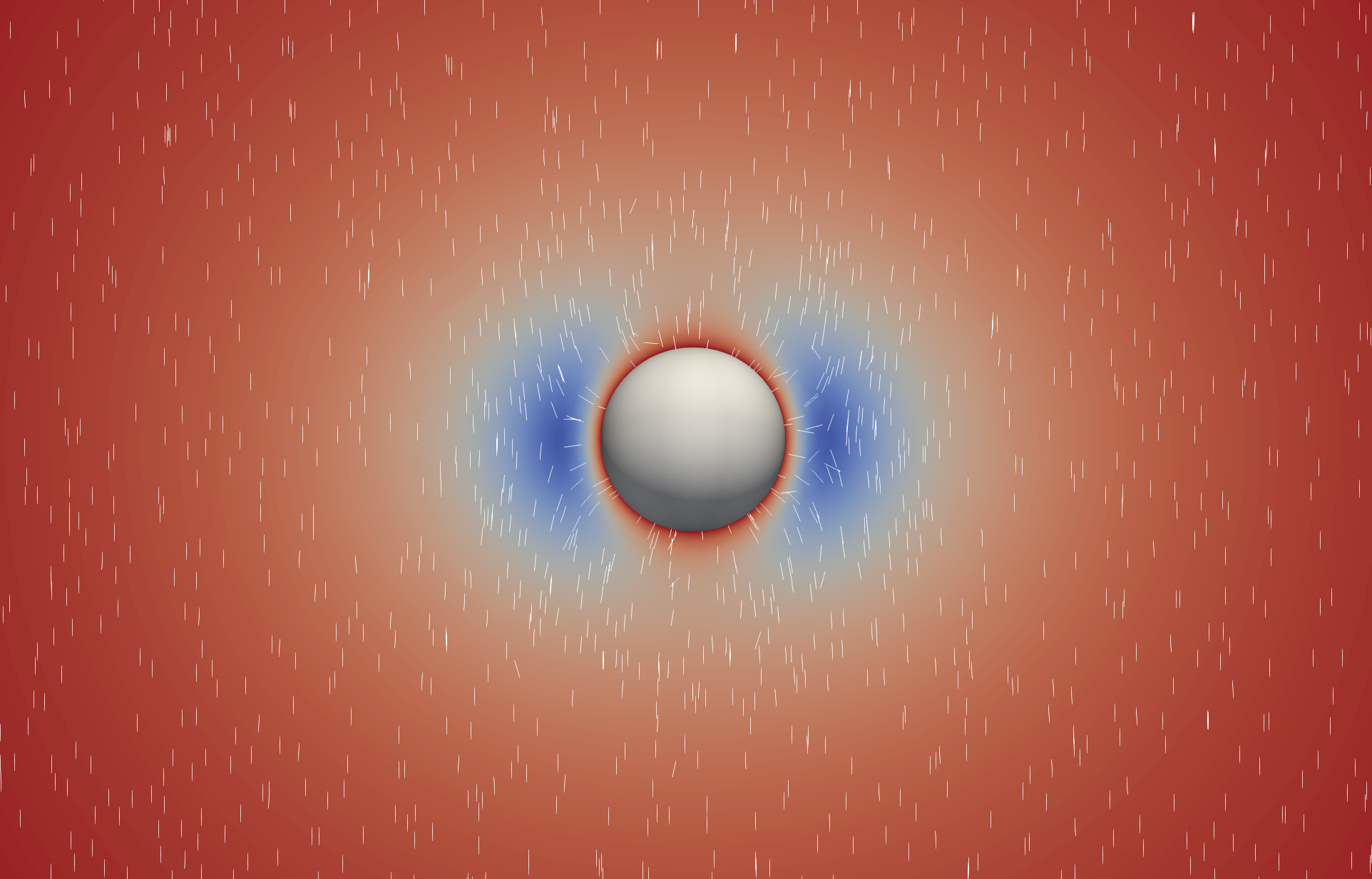} 
\includegraphics[scale=0.15]{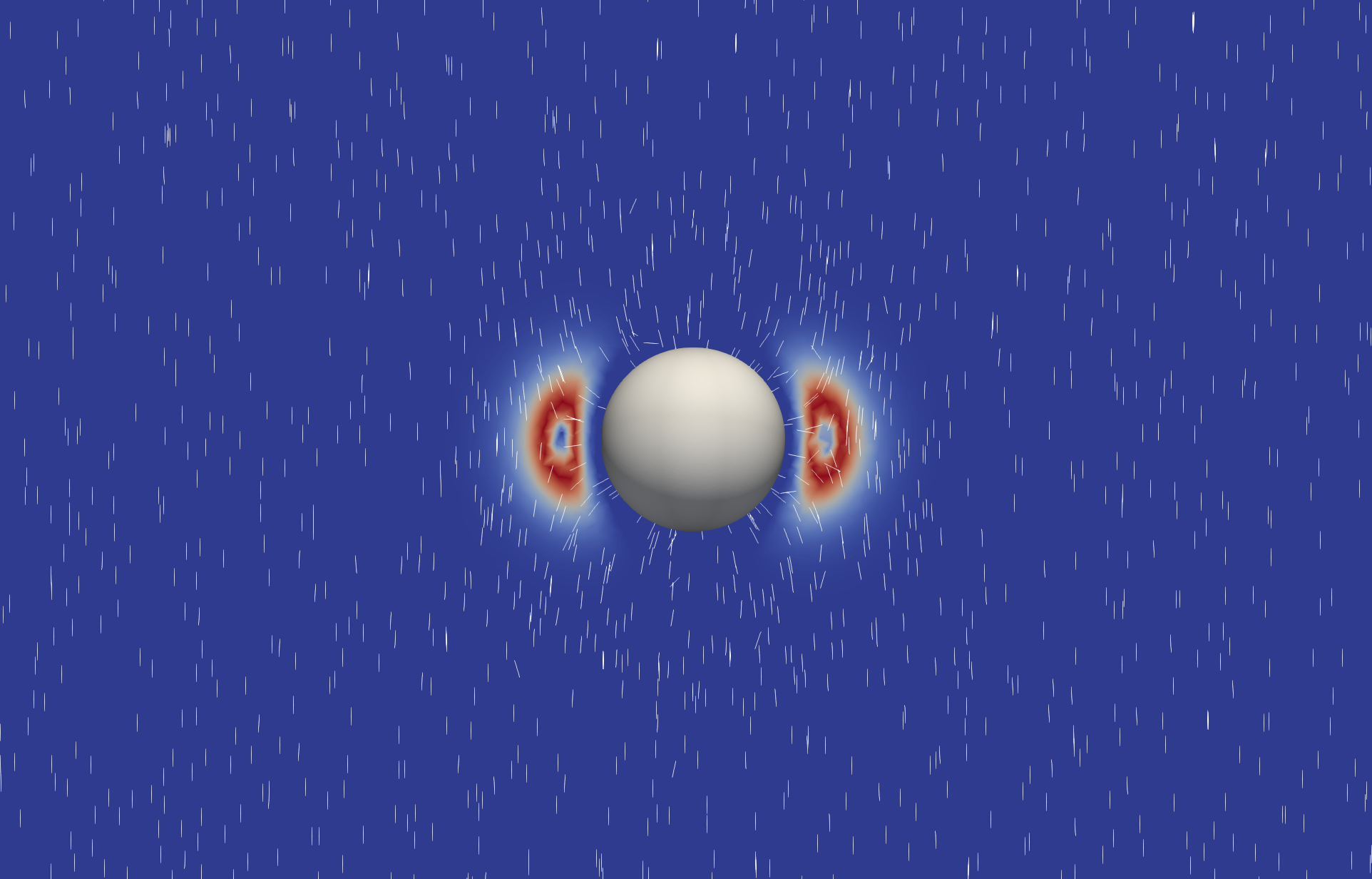} 
\includegraphics[scale=0.15]{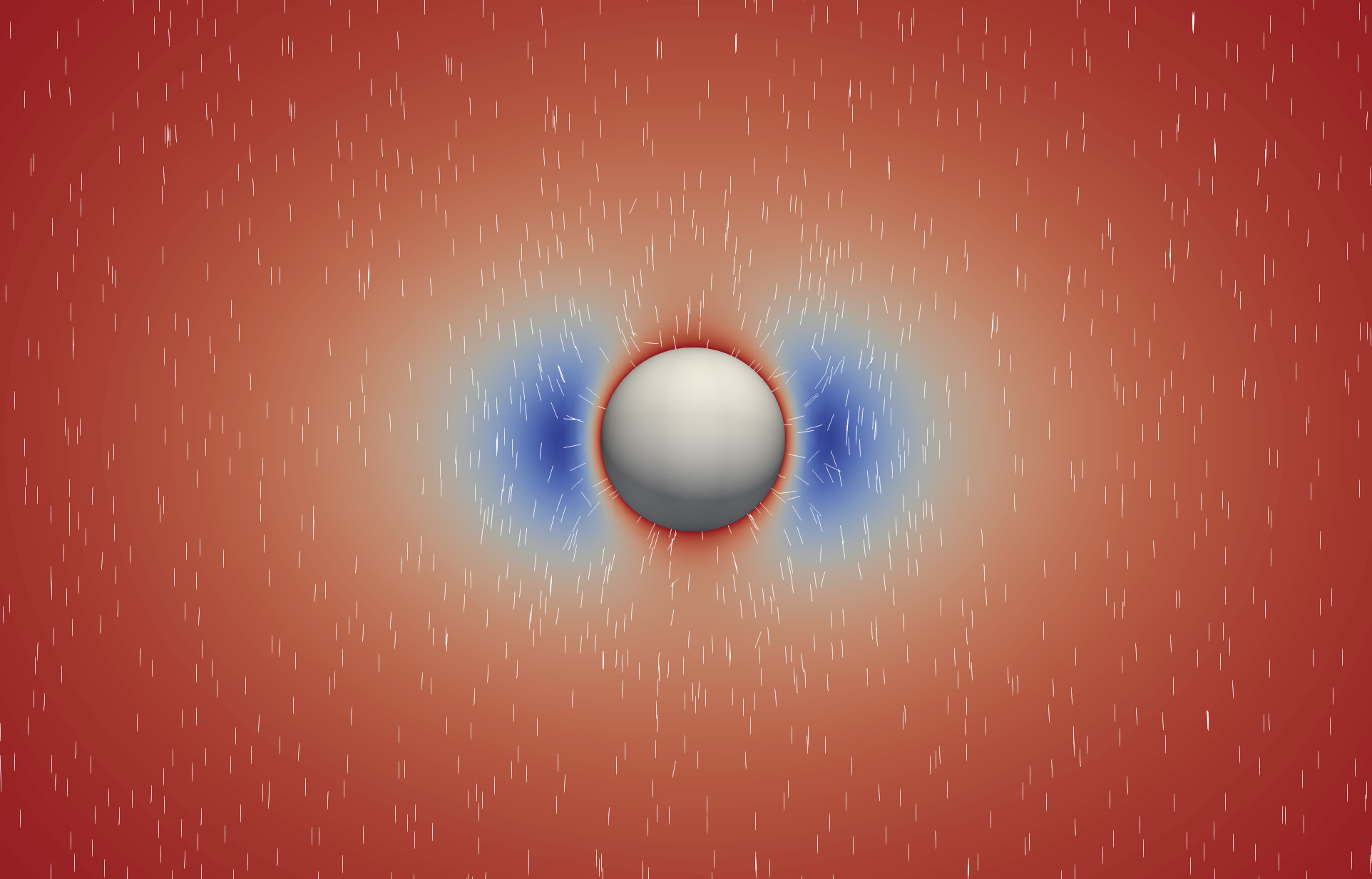} 
\includegraphics[scale=0.15]{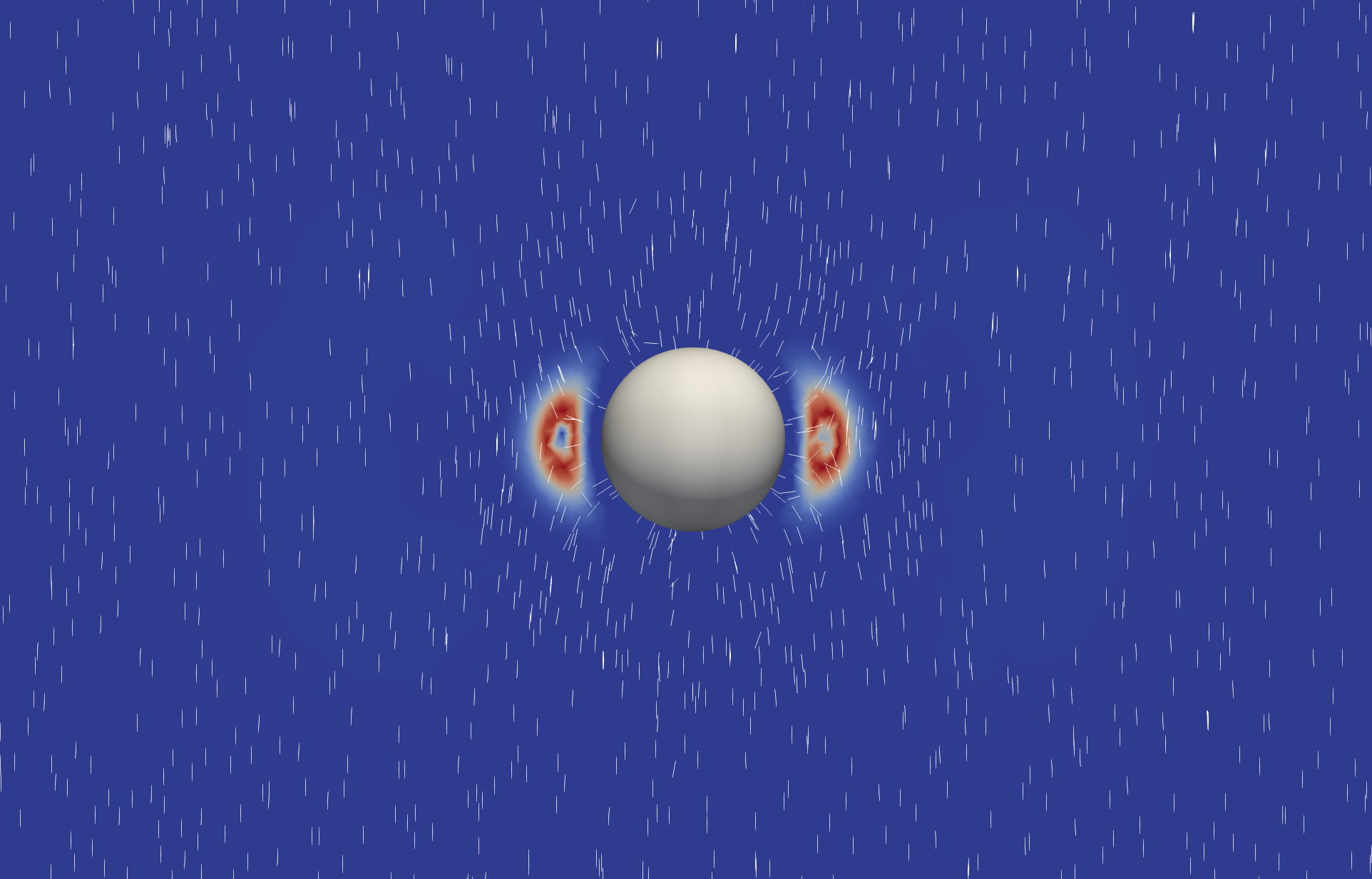} 
\caption{Plots of the norm $|Q|$ (left column) and biaxiality parameter $\beta(Q)$ (right column) for $k=0,-0.5,-0.99$ (top to bottom). 
The lines indicate the director field (dominating eigenvector) of $Q$.
The color coding is identical for all four norm images and goes from $0.18$ (dark blue) to $0.82$ (red).
The biaxiality parameter varies between $0$ (uniaxial, dark blue) and $1$ (maximally biaxial, red).
}
\label{fig:norm_biax_k_neg_0_05_099}
\end{center}
\end{figure}

Of particular interest is the defect region of a solution $Q$. 
There are different ways to characterize these regions.
One approach uses the \emph{biaxiality parameter} $\beta(Q)$ which is defined as
\begin{align*}
    \beta(Q)
    \ := \
    1 - 6 \frac{\tr(Q^3)^2}{\tr(Q^2)^3}
    \, .
\end{align*}
We remark that $\beta(Q)\in [0,1]$ with $\beta(Q)=0$ if and only if $Q$ is uniaxial \cite[Lemma 13]{MaZar}.
It can therefore be seen as a measure of the deviation of $Q$ from being uniaxial and provides a useful indication where the defects are located, see Figure \ref{fig:norm_biax_k_0_5_10_20} and \ref{fig:norm_biax_k_neg_0_05_099}.
Another possibility to detect defects is by explicitly computing the director field associated to $Q$.
An abrupt change of direction indicates a defect, see again Figure \ref{fig:norm_biax_k_0_5_10_20} and \ref{fig:norm_biax_k_neg_0_05_099}.

\begin{figure}[H]
\begin{center}
\includegraphics[scale=0.15]{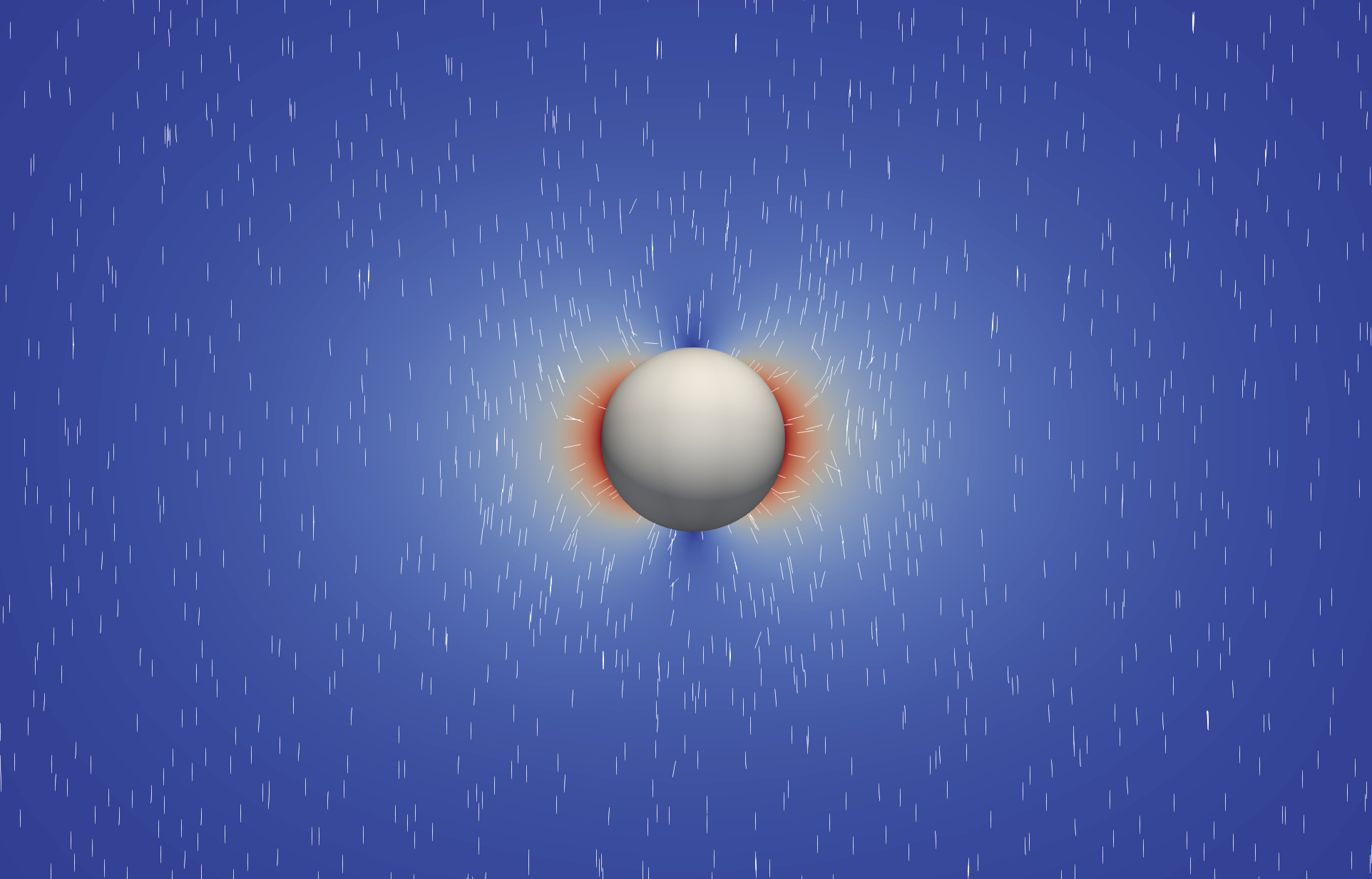} 
\includegraphics[scale=0.15]{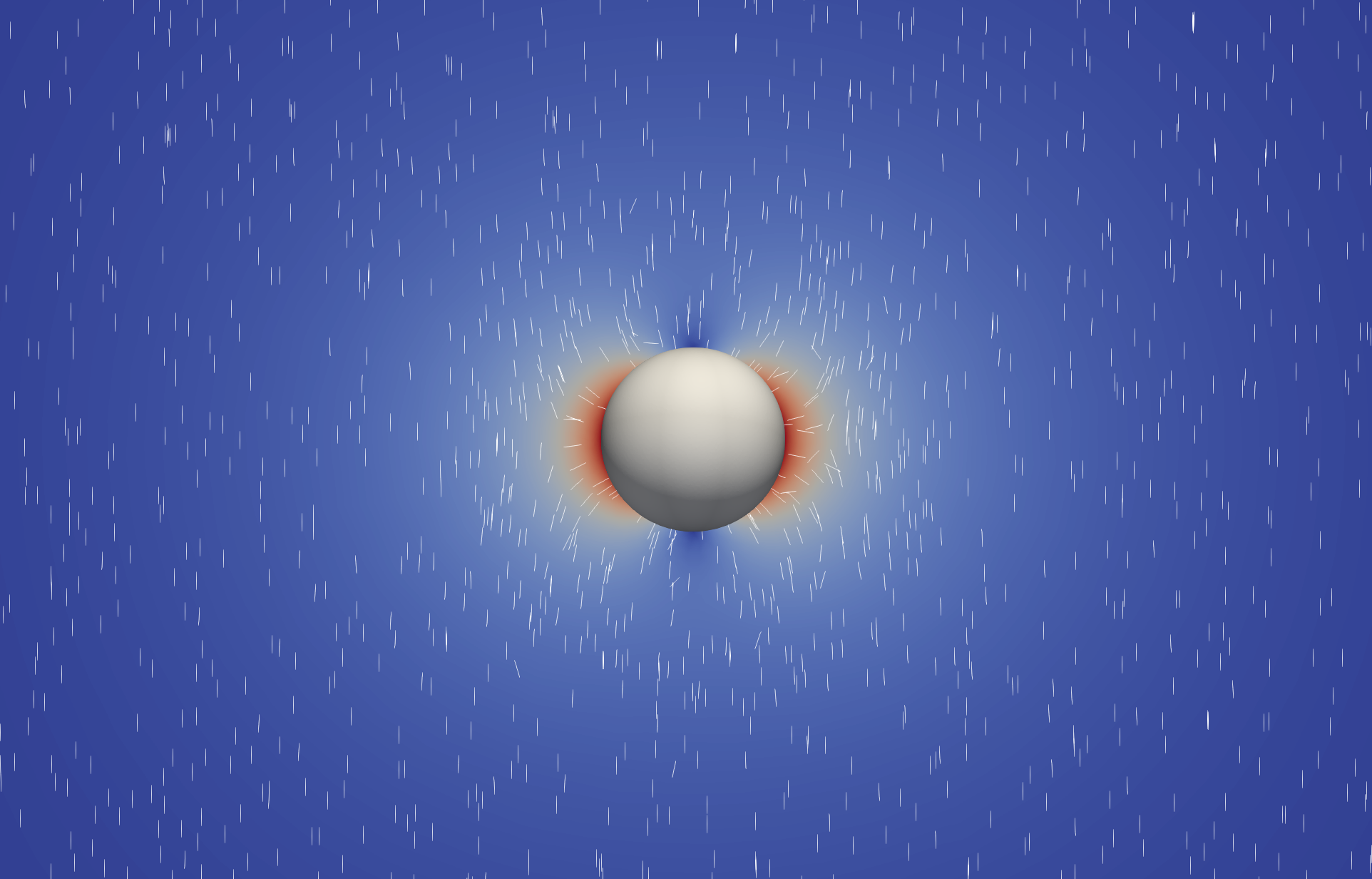} 
\includegraphics[scale=0.15]{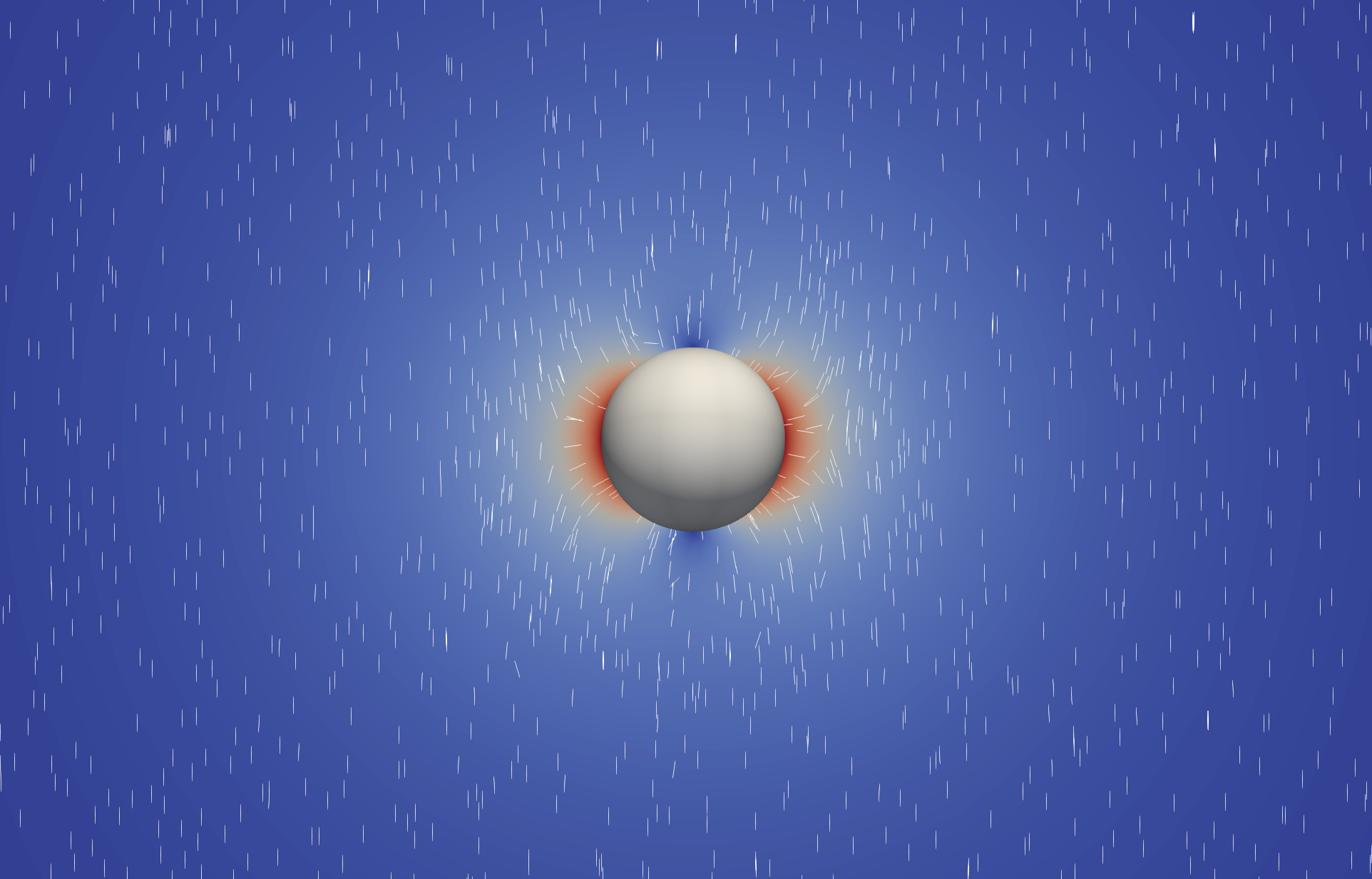} 
\includegraphics[scale=0.15]{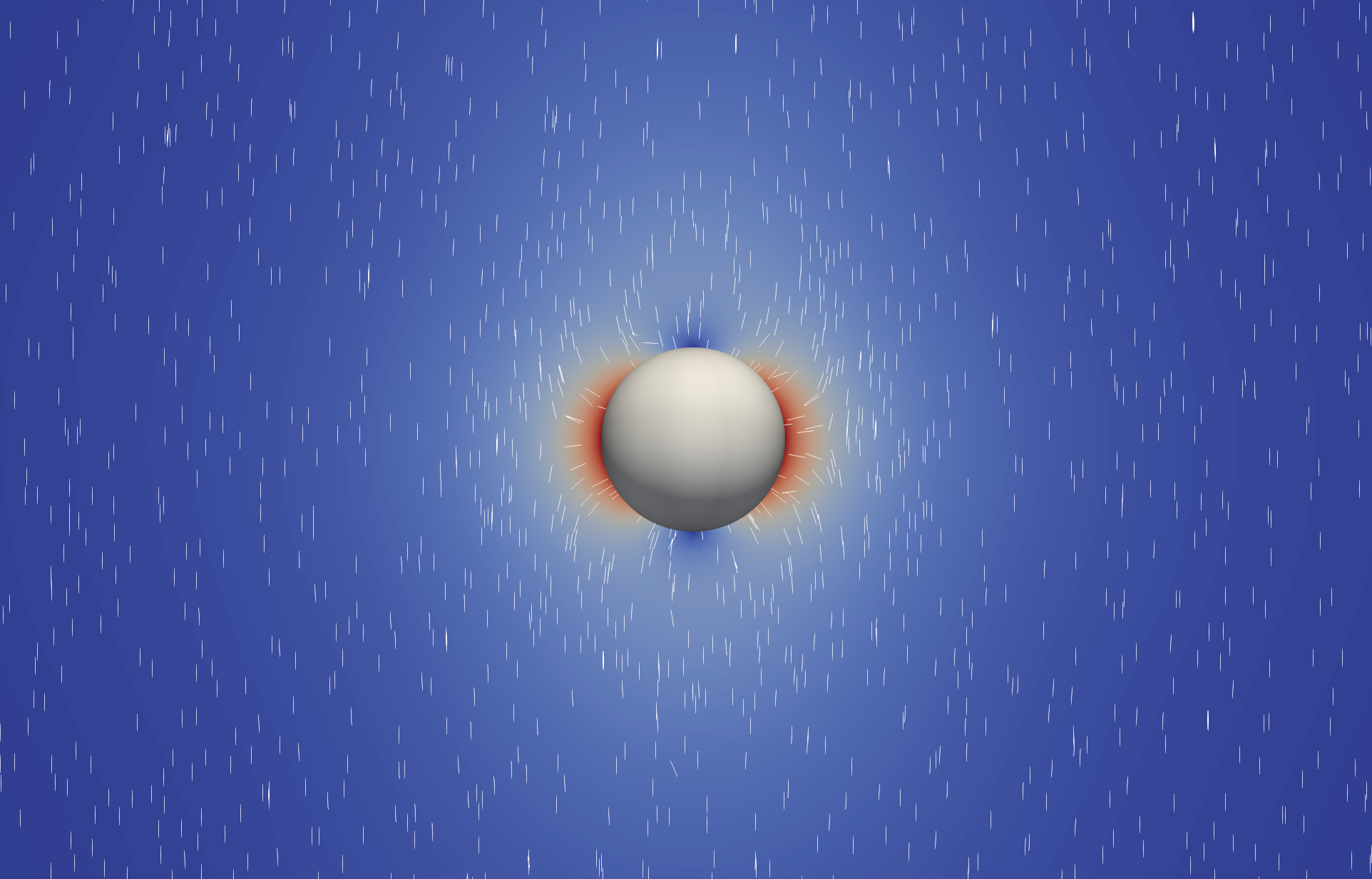} 
\includegraphics[scale=0.15]{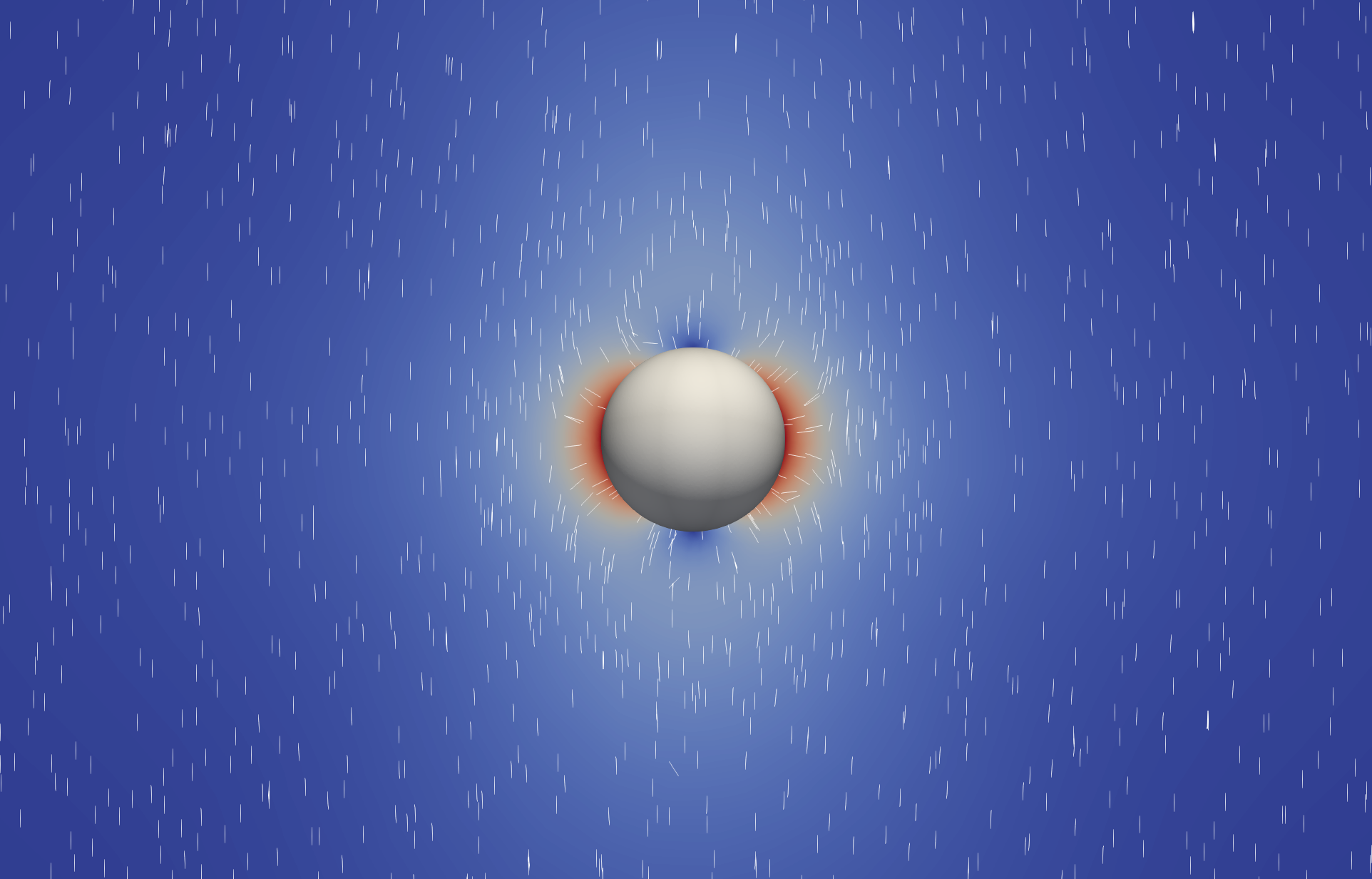} 
\includegraphics[scale=0.15]{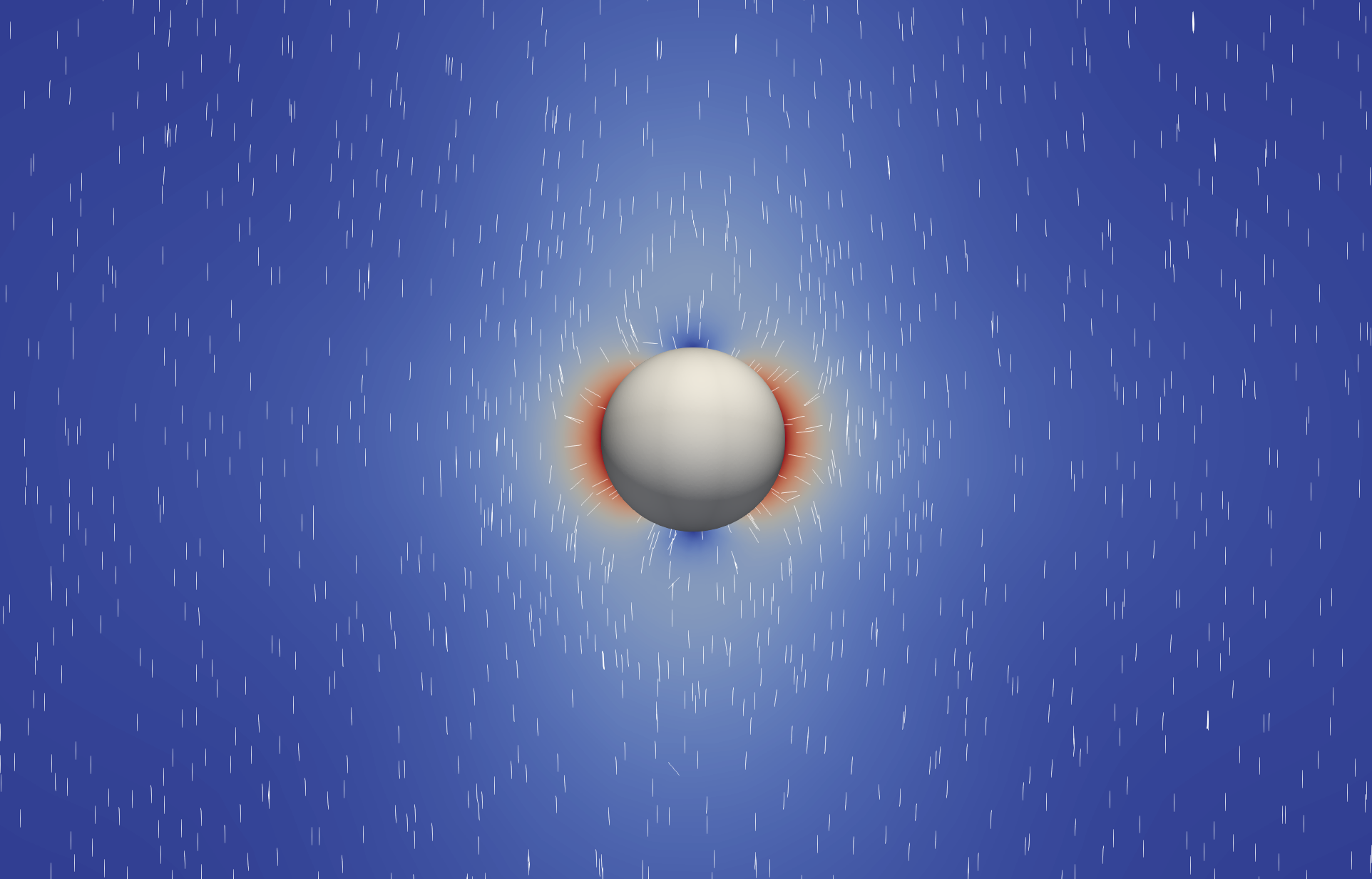} 
\caption{Plots of the distance of $Q$ from $Q_\infty$ given by $|Q-Q_\infty|$ for $k=-0.99,-0.5,0,5,10,20$ (left to right, top to bottom). 
The lines indicate the director field (dominating eigenvector) of $Q$.
The color coding is identical for all six images and goes from $0$ (dark blue) to $1.41$ (red).
}
\label{fig:distQinf_k_m1_0_5_10_20}
\end{center}
\end{figure}

\begin{figure}[H]
\begin{center}
\includegraphics[scale=0.20]{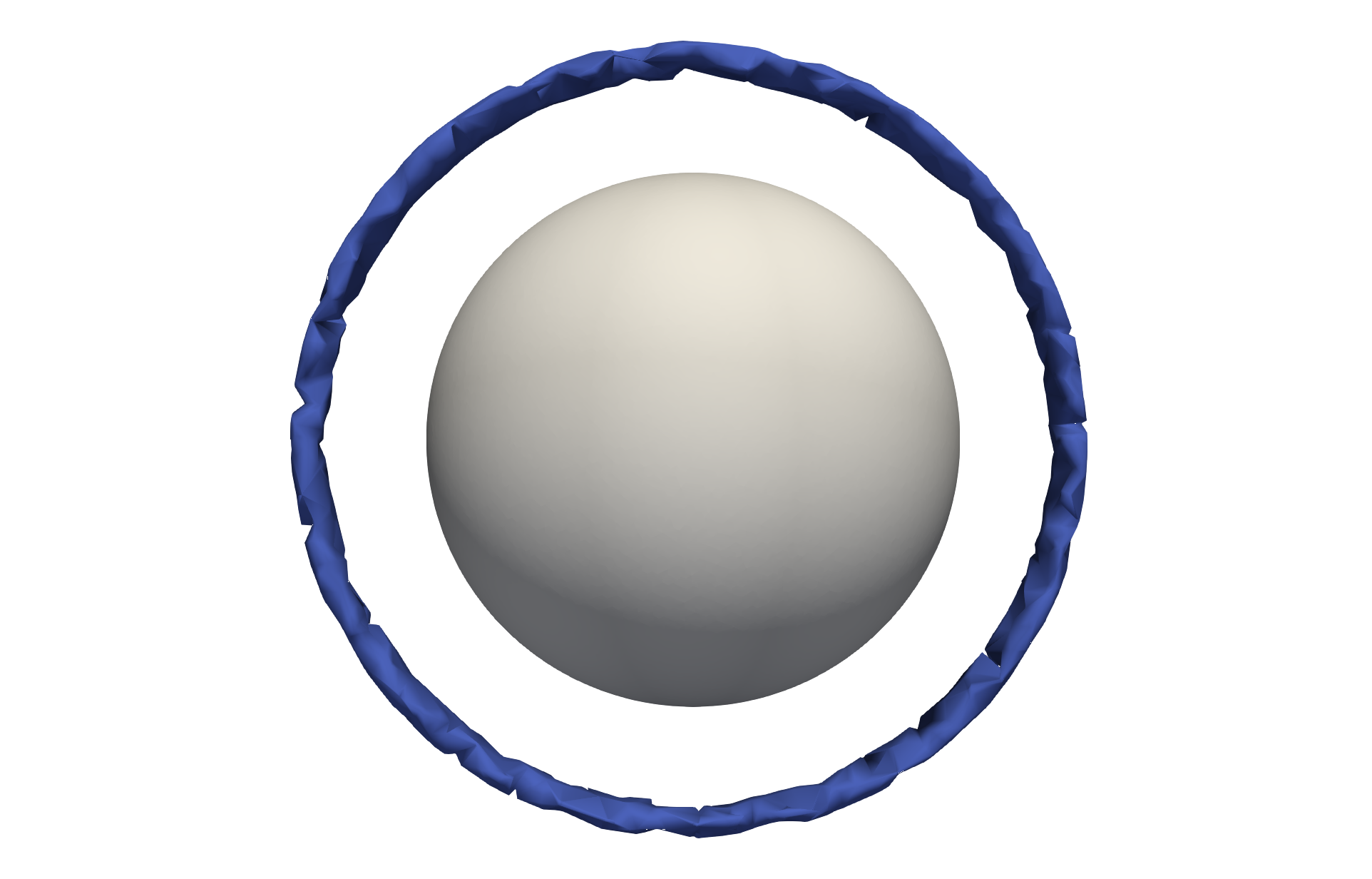} 
\includegraphics[scale=0.20]{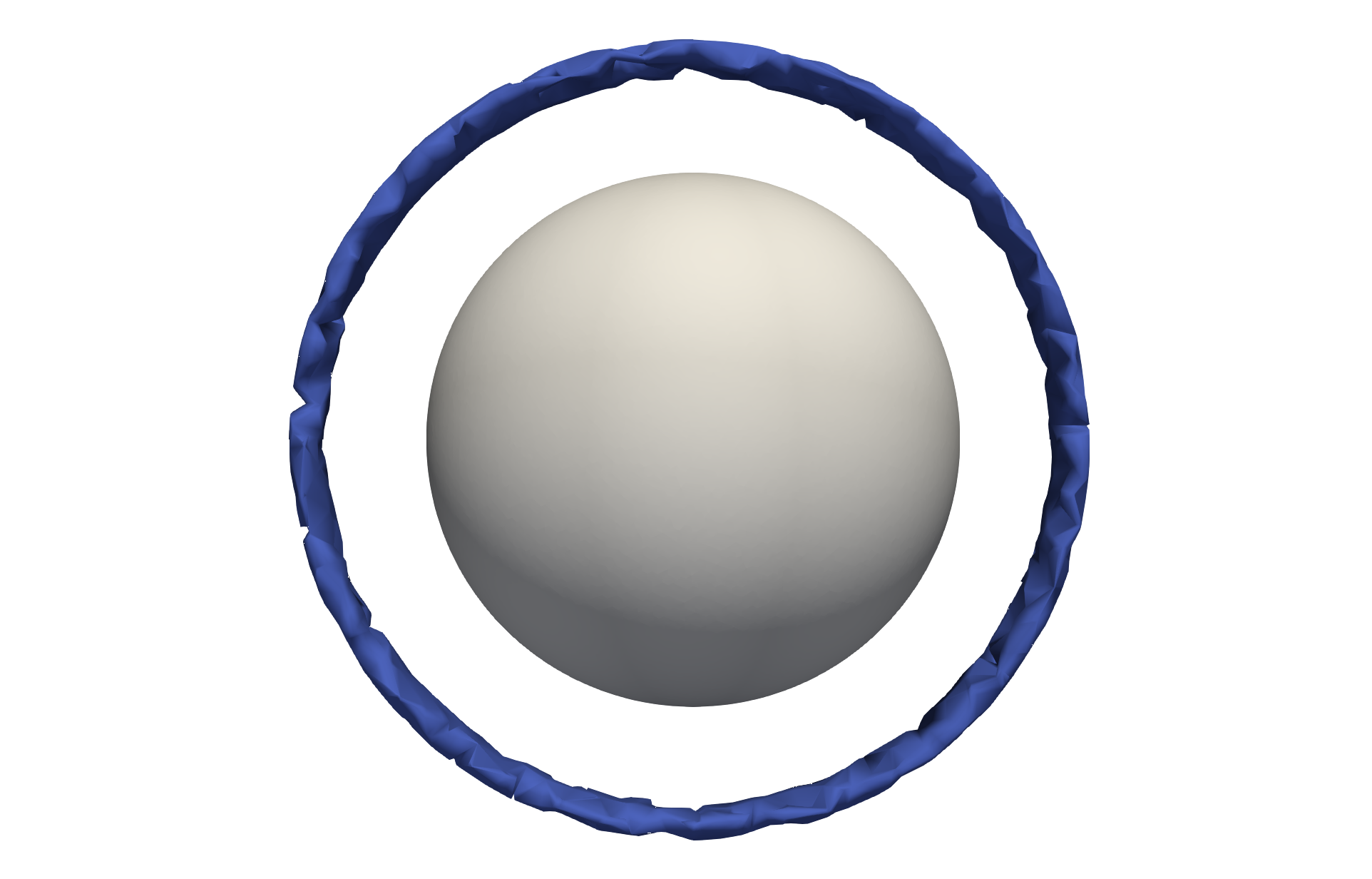} 
\includegraphics[scale=0.20]{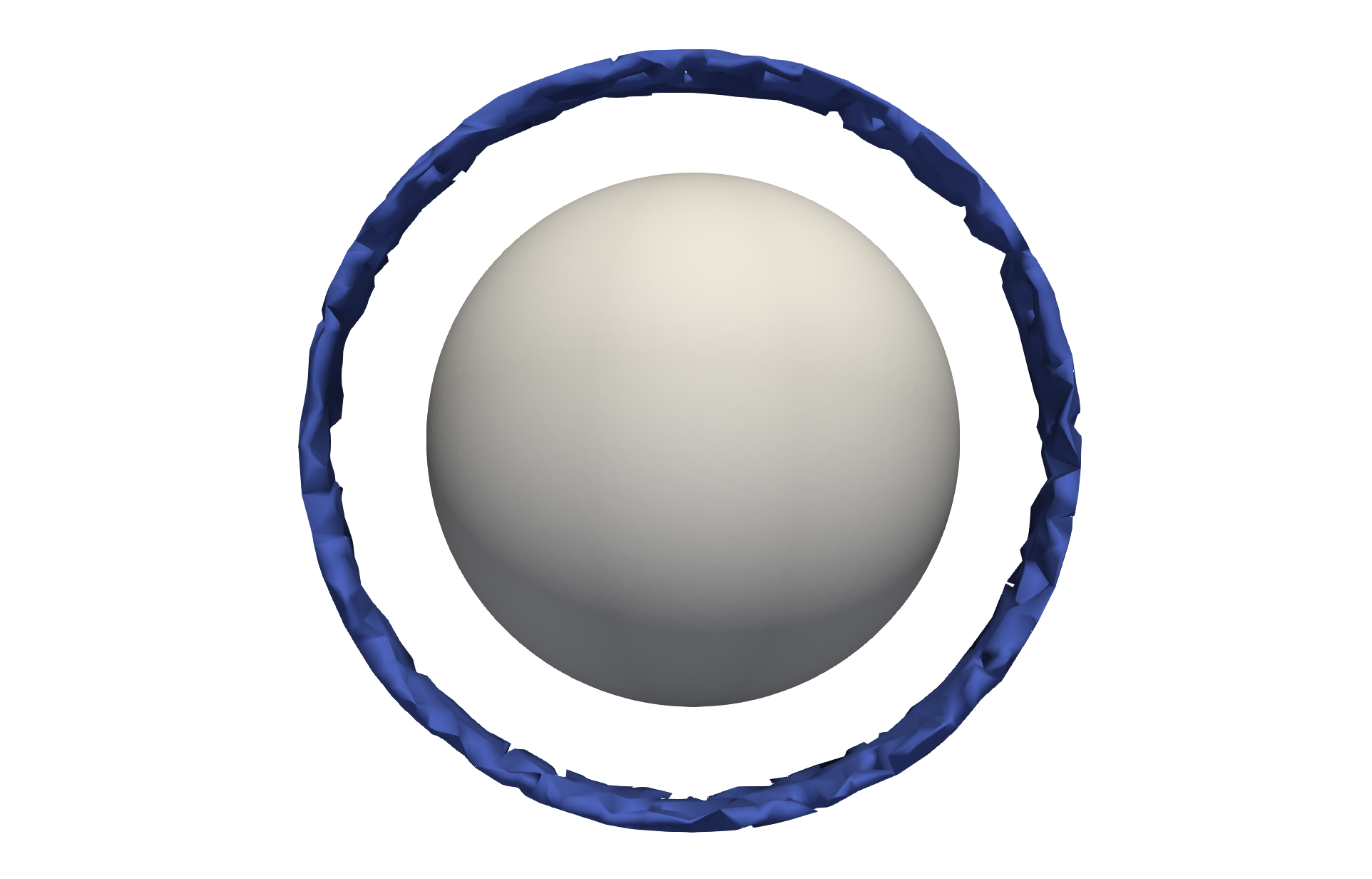} 
\caption{Plots of the ``Saturn ring'' surrounding the spherical colloid for $k=-0.99,0,20$ (top to bottom). 
The ring shown corresponds to a difference of $0.08$ in the leading eigenvalues of $Q$.
}
\label{fig:SR_k_m1_0_20}
\end{center}
\end{figure}

For all tested values of $k$, we find a Saturn ring configuration around the colloid.
The ring defect (surrounded by a biaxial region) stays at the equatorial plane and even the radius stays roughly constant (only decreases marginally with $k$), see Figure \ref{fig:SR_k_m1_0_20}. 
However, the overall size of the region in which $Q$ is biaxial increases with $k$, compare Figure \ref{fig:norm_biax_k_0_5_10_20}. 

One can observe that for $k=0$ the norm is smallest in the defect region, here $Q$ is closest to isotropic.
This effect becomes even stronger for negative $k$ (Figure \ref{fig:norm_biax_k_neg_0_05_099}).
With increasing $k$, the norm inside the Saturn ring increases, while the norm above (resp.\ below) the North pole (resp.\ South pole) decreases, as observed in Figure \ref{fig:norm_biax_k_0_5_10_20}. 
The region in which the norm of $Q$ is small forms an ellipsoidal shape that changes from oblate (for small $k$) to prolate (for larger $k$).

As depicted in Figure \ref{fig:distQinf_k_m1_0_5_10_20}, the decay of $|Q-Q_\infty|$ is not isotropic, but varies as a function of the angle between $e_3$ and the normal to the colloid $\nu$:
For small angles (i.e.\ vertical directions) we observe a rapid decay for small $k$ that becomes slower for larger $k$.
Large angles (i.e.\ horizontal directions) present the opposite phenomenon, that is, the decay increases with increasing $k$.

\bibliographystyle{acm}
\bibliography{bibliography}

\end{document}